\documentclass[3pt]{elsarticle}
\pdfoutput=1
\usepackage{amsthm}
\usepackage[intlimits]{amsmath}
\usepackage{amssymb}
\usepackage{amsthm}
\usepackage{amsfonts}
\usepackage{bm}
\usepackage{mathdots}
\usepackage[normalem]{ulem}

\usepackage[utf8]{inputenc} 
\usepackage[T1]{fontenc}
\usepackage{times}
\usepackage{mathptmx}  
\DeclareSymbolFont{cmletters}{OML}{cmr}{m}{n}
\DeclareMathAlphabet{\mathcal}{OMS}{cmsy}{m}{n} 

\journal{Advances in Engineering Software}
\ifdefined\copyReadyOn
\usepackage{ecrc}
\volume{00}
\firstpage{1}
\runauth{G. Beer}
\journalname{\journal}
\jid{CG}
\jnltitlelogo{CG}
\CopyrightLine{2015}{Published by Elsevier Ltd.}
\fi

\usepackage{tikz}
\usepackage{tikz-3dplot}
\usepackage{pgfplots}
\usepackage{pgfplotstable}
\usepackage{tikzscale}
\usepackage{standalone} 
\usepackage{hf-tikz}

\usepackage{forloop}
\usepackage{arrayjobx}

\usepackage{algorithmic}
\usepackage{algorithm}

\usepackage{booktabs} 
\usepackage{calc}
\usepackage{xcolor,pdfcomment,soul} 
\usepackage{xspace} 
\usepackage{multirow}

\usepackage{scrlayer-scrpage}
\usepackage{setspace}

\usepackage[labelformat=simple]{subcaption}   	
\usepackage{pdfpages} 

\usepackage{tabto} 

\usepackage{makeidx}
\makeindex

\usepackage{colortbl}

\usepackage{enumerate}

\usepackage{overpic}

\usepackage{listings} 

\usepackage[intlimits]{amsmath}
\usepackage{amssymb} 
\usepackage{amsthm}
\usepackage{amsfonts,mathrsfs}
\usepackage{bm}
\usepackage{mathtools} 

\usepackage{siunitx} 

\usepackage{xfrac}
\usepackage{listings}
\usepackage{color}
\definecolor{dkgreen}{rgb}{0,0.6,0}
\definecolor{gray}{rgb}{0.5,0.5,0.5}
\definecolor{mauve}{rgb}{0.58,0,0.82}

\usepackage[utf8]{inputenc} 
\usepackage[german,english]{babel}  
\usepackage[T1]{fontenc} 
\usepackage{lmodern}

\usepackage{geometry}
\usepackage{overpic}
\graphicspath{{pics/}{tmp/}}

\newboolean{inGray}
\setboolean{inGray}{false}

\newcommand{\addtoindex}[2][]{
    \ifthenelse { \equal{#1}{} }
    {#2\index{#2}\xspace}%
    {#2\index{#1}\xspace}%
}

\newcommand{\myVec}[1]{\mathbf{#1}}

\newcommand\uu{\xi} 			
\newcommand\vv{\eta} 			

\newcommand\abscissa{\tilde{\uu}}			
\providecommand\url[1]{\emph{#1}}

\newcommand\pt[1]{\boldsymbol{#1}}

\newenvironment{mytable}[4]%
{
  \begin{table}[#1]
    \def\mycap{#2}
    \def\mylabel{#3}
    \centering
    \begin{tabular}{#4}
      \toprule
}
{
  \bottomrule
  \end{tabular}
  \caption{\mycap}
  \label{\mylabel}
  \end{table}
}
\newcommand{\mytableheader}[1]{
       #1 \\ \midrule
}

\newenvironment{myalgorithm}[2]%
{
  \begin{algorithm}
    \caption{#1}
    \label{#2}
    \begin{algorithmic}[1]
}   
{
\end{algorithmic}
\end{algorithm}
}

\newtheoremstyle{myremark}
{3pt}
{3pt}
{}
{}
{\itshape}
{:}
{.5em}
{}
\theoremstyle{myremark}
\newtheorem*{remark}{Remark}

\newcommand{\myeqref}[1]{equation~(\ref{#1})}	
\newcommand{\myfigref}[1]{Figure~\ref{#1}}

\newcounter{footnoteNumber} 

\DeclareCaptionFormat{listing}{#1#2#3}
\newcommand{\myalignatsinglelabel}[3]%
{%
    \begin{equation}
        \label{#2}
        \begin{alignedat}{#1} 
            #3
        \end{alignedat}
    \end{equation}
}

\newcommand{\myalignat}[2]%
{%
  \begin{alignat}{#1} 
    #2
   \end{alignat}
}

\usepackage{pgf,pgffor,pgfmath}
\pgfplotsset{compat=1.10}

\usetikzlibrary{external}
\usetikzlibrary{positioning,calc,intersections,arrows,3d}
\usetikzlibrary{shapes,snakes,fit,spy}
\usepgfplotslibrary{patchplots}
\usetikzlibrary{matrix}

\usepgfplotslibrary{groupplots}
\usepgfplotslibrary{fillbetween}

\def\pgfplotfontsizetitle{\small}
\def\pgfplotfontsizelegend{\small}
\def\pgfplotfontsize{\small}
\def\pgfplotfontsizetiny{\scriptsize}

\def\tikzfontsizetiny{\scriptsize}

\pgfplotsset{
  mystyle/.style ={%
    grid = major,
    every tick label/.append style={font=\pgfplotfontsizetiny},
    every axis label/.append style={font=\pgfplotfontsize},
    legend style={font=\pgfplotfontsizelegend},
    label style={font=\pgfplotfontsize},
    title style={font=\pgfplotfontsizetitle},
    /pgf/number format/set thousands separator = {}, 
  }
}%

\makeatletter
\pgfplotsset{
    myIgnoreRowModulo2/.style args={#1}{%
        /pgfplots/x filter/.code={%
        \let\xValue\pgfmathresult 
        \pgfmathparse{int(mod(int(\coordindex),int(2))} \pgfmathresult 
        \ifnum#1=\pgfmathresult
            \def\pgfmathresult{} 
        \else
            \pgfmathparse{\xValue} \pgfmathresult 
        \fi
        }
    } 
}
\makeatother

\colorlet{drawblue}      {blue!80!white}
\colorlet{drawred}       {red!80!white}
\colorlet{drawgray}      {gray}
\definecolor{drawgreen}  {RGB}{44,162,95}
\colorlet{drawpurple}    {purple}
\colorlet{draworange}    {orange}
\colorlet{drawlime}      {lime!80!black}
\colorlet{drawartichoke} {yellow!60!black}
\colorlet{TUGgray}{black!15}
\definecolor{TUGred}{RGB}{247,1,70}
\definecolor{IFBblue}{RGB}{51,112,169}

\definecolor{basisColor1}{RGB}{59,76,192}
\definecolor{basisColor2}{RGB}{87,117,225}
\definecolor{basisColor3}{RGB}{119,154,247}
\definecolor{basisColor4}{RGB}{152,185,255}
\definecolor{basisColor5}{RGB}{184,208,249}
\definecolor{basisColor6}{RGB}{195,209,230}	
\definecolor{basisColor7}{RGB}{247,200,190}	
\definecolor{basisColor8}{RGB}{247,187,160}
\definecolor{basisColor9}{RGB}{244,154,123}
\definecolor{basisColor10}{RGB}{229,112,88}
\definecolor{basisColor11}{RGB}{203,62,56}
\definecolor{basisColor12}{RGB}{180,4,38}

\definecolor{basisColor8sw}{RGB} {189,189,189}
\definecolor{basisColor9sw}{RGB} {150,150,150}
\definecolor{basisColor10sw}{RGB}{115,115,115}
\definecolor{basisColor11sw}{RGB}{91,91,91}
\definecolor{basisColor12sw}{RGB}{37,37,37}

\colorlet{myblue}    {blue}
\colorlet{myred}     {red}
\colorlet{mygreen}   {drawgreen}
\colorlet{mypurple}  {purple}
\colorlet{myorange}  {orange}

\ifthenelse{\boolean{inGray}}{ 

\pgfplotscreateplotcyclelist{mycyclelist}{
  basisColor12sw,semithick, every mark/.append style={solid, fill=.},mark=*           \\%
  basisColor11sw,semithick, every mark/.append style={solid, fill=.},mark=square*     \\%
  basisColor10sw,semithick, every mark/.append style={solid, fill=.},mark=triangle*   \\%
  basisColor9sw , semithick,every mark/.append style={solid, fill=.},mark=diamond*    \\%
  basisColor8sw , semithick,every mark/.append style={solid, fill=.},mark=otimes    \\%
}
}{

\pgfplotscreateplotcyclelist{mycyclelist}{
  basisColor12,semithick,   every mark/.append style={solid, fill=.},mark=*           \\%
  basisColor11,semithick,   every mark/.append style={solid, fill=.},mark=square*     \\%
  basisColor10,semithick,   every mark/.append style={solid, fill=.},mark=triangle*   \\%
  basisColor9, semithick,   every mark/.append style={solid, fill=.},mark=diamond*    \\%
  basisColor8, semithick,   every mark/.append style={solid, fill=.},mark=otimes    \\%
}
}

\ifthenelse{\boolean{inGray}}{ 
\pgfplotscreateplotcyclelist{mycyclelistcompare}{
  basisColor12sw, semithick, loosely dashdotdotted     \\%
  basisColor11sw, semithick, loosely dashed      		\\%
  basisColor10sw, semithick, dashed \\
  basisColor9sw,  semithick, solid               	\\%
  basisColor12sw, only marks,   every mark/.append style={solid, fill=.},mark=otimes    \\%
  basisColor11sw, only marks,   every mark/.append style={solid, fill=.},mark=triangle*   \\%
  basisColor10sw, only marks,   every mark/.append style={solid, fill=.},mark=square*     \\%
  basisColor9sw,  only marks,   every mark/.append style={solid, fill=.},mark=*          \\%
}

}{
\pgfplotscreateplotcyclelist{mycyclelistcompare}{
  basisColor12, semithick, loosely dashdotdotted     \\%
  basisColor11, semithick, loosely dashed      		\\%
  basisColor10, semithick, dashed \\
  basisColor9,  semithick, solid               	\\%
  basisColor12, only marks,   every mark/.append style={solid, fill=.},mark=otimes    \\%
  basisColor11, only marks,   every mark/.append style={solid, fill=.},mark=triangle*   \\%
  basisColor10, only marks,   every mark/.append style={solid, fill=.},mark=square*     \\%
  basisColor9,  only marks,   every mark/.append style={solid, fill=.},mark=*          \\%
}
}

\ifthenelse{\boolean{inGray}}{ 
\tikzset{mycyclelistcompareReferenceA/.style={basisColor12sw,solid}}
\tikzset{mycyclelistcompareTestA/.style={basisColor12sw,only marks,mark=otimes}}
}{
\tikzset{mycyclelistcompareReferenceA/.style={basisColor12,solid}}
\tikzset{mycyclelistcompareTestA/.style={basisColor12,only marks,mark=otimes}}
}

\tikzset{helpline/.style={thin,dashed}}
\tikzset{labelline/.style={thin}}
\tikzset{referencePath/.style={dotted,very thick}}
\tikzset{helparrow/.style={thin,arrows={-latex}}}
\tikzset{axis/.style={thin,arrows={->}}}
\tikzset{force/.style={thick,arrows={->}}}
\tikzset{forceInverse/.style={thick,arrows={<-}}}
\tikzset{Gamma/.style={ultra thick}}
\tikzset{controlPoly/.style={draw=black}}
\tikzset{GammaFill/.style={fill=lightgray,fill opacity=0.5}}
\tikzset{colorDiri/.style={drawgreen}}
\tikzset{GammaFillDiri/.style={fill=drawgreen,fill opacity=0.5}}
\tikzset{initialgrid/.style={thin,gray}}
\tikzset{addgridline/.style={dashed,gray}}
\tikzset{trimmingcurve/.style={thick}}
\tikzset{boundingbox/.style={thick, dotted}}
\tikzset{parameterSpace/.style={ }}
\tikzset{basisfunction/.style={very thick,smooth}}
\tikzset{bspline/.style={very thick,smooth,red}}
\tikzset{intersectioncurve/.style={dashed,thick}}
\tikzset{integrationRegionEdge/.style={dashed}}
\tikzset{pointer/.style={arrows={-latex}}}

\tikzstyle{anode}= [circle, inner sep=1.3pt, draw, fill=black]
\tikzstyle{gausspoint}=[shape=circle,draw=black,fill=black,inner sep=1.1pt]
\tikzstyle{controlPoint}=[shape=circle,draw=black,fill=black,thin,inner sep=0pt,minimum size=1.5mm]
\tikzstyle{abscissaPoint}=[shape=circle,draw=black,fill=white,thin,inner sep=0pt,minimum size=1.5mm]
\tikzstyle{anchorPoint}=[shape=circle,draw=black,fill=black,thin,inner sep=0pt,minimum size=1.5mm]
\tikzstyle{anchorPointDeg}=[shape=circle,draw=black,fill=TUGred,thin,inner sep=0pt,minimum size=1.5mm]
\tikzstyle{anchorPointDegD}=[shape=cross out,thick,draw=black,inner sep=0pt,minimum size=1.5mm]
\tikzstyle{trimmingIntersectionPoint}=[shape=cross out,thick,draw=black,inner sep=0pt,minimum size=1.5mm]

\tikzset{%
  highlight/.style={rectangle,rounded corners,fill=red!60,draw,fill opacity=0.125,thick,inner sep=0pt}
}

\usepackage{colortbl}

\def\trianglecolor{black}
\newcommand{\upperSlopeTriangle}[4] 	
{
	\def\trianglecolor{black}
	\addplot[forget plot, domain=#3:#4,color=\trianglecolor,samples=2]{  #2 / (x^#1) } node (A1) [pos=1] {}; 
	\addplot[forget plot, domain=#3:#4,color=\trianglecolor,samples=2]{   #2  / (#3^#1)} node (A2) [pos=1] {} node [anchor=south,pos=0.5,black] {\tikzfontsizetiny $1$};
	\draw[color=\trianglecolor] (A1.center) -- (A2.center) node [anchor=west,pos=0.5,black] {\tikzfontsizetiny #1};
}

\newcommand{\lowerSlopeTriangle}[4] 	
{
	\def\trianglecolor{black}
	\addplot[forget plot, domain=#3:#4,color=\trianglecolor,samples=2]{  #2 / (x^#1) } node (A1) [pos=0] {}; 
	\addplot[forget plot, domain=#3:#4,color=\trianglecolor,samples=2]{   #2  / (#4^#1)} node (A2) [pos=0] {} node [anchor=north,pos=0.5,black] {\tikzfontsizetiny $1$};
	\draw[color=\trianglecolor] (A1.center) -- (A2.center) node [anchor=east,pos=0.5,black] {\tikzfontsizetiny #1};
}

\newcommand{\myaddgraphic}[5]
{
 \node[anchor=south west,inner sep=0] (image) {\phantom{\includegraphics[#2]{#1}}};
  \begin{scope}[x={(image.south east)},y={(image.north west)}]
      
      \begin{scope}
          
          #5
          
          \node[anchor=south west,inner sep=0] {\includegraphics[#2]{#1}};
      \end{scope} 
      
      #4
      
      \pgfmathparse{int(#3)} \let\gridIndicator\pgfmathresult
      \ifthenelse{ \gridIndicator = 1 }
      {
          \draw[help lines,xstep=.1,ystep=.1] (0,0) grid (1.001,1.001);
          \foreach \x in {1,...,9} { \node [anchor=north] at (\x/10,0) {\x};}
          \foreach \y in {1,...,9} { \node [anchor=east] at (0,\y/10) {\y};}
      }{}
      
  \end{scope}    
}

\tikzstyle{reverseclip}=[insert path={(current page.north east) --
    (current page.south east) --
    (current page.south west) --
    (current page.north west) --
    (current page.north east)}
]

\newcounter{itR}
\newcommand{\bsplinevalue}[5] 
{                    				
	\newarray\vKnots
	\newarray\vN
	\newarray\vNumeratorL
	\newarray\vNumeratorR
	\newarray\vSave

	\readarray{vKnots}{#1}
	\readarray{vN}{1}
	\readarray{vSave}{0}
	\readarray{vNumeratorL}{0}
	\readarray{vNumeratorR}{0}

	\foreach \j in {1,...,#2}
	{        
		\pgfmathparse{ int(#4+\j+1) } \checkvKnots(\pgfmathresult)
		\pgfmathsetmacro{\numR}{\cachedata-#3}
		
		\pgfmathparse{ int(#4-\j+2) } \checkvKnots(\pgfmathresult)
		\pgfmathsetmacro{\numL}{#3-\cachedata} 					

		\expandarrayelementtrue
		\pgfmathparse{ int(\j+1) }
		\vNumeratorL(\pgfmathresult)={\numL}
		\vNumeratorR(\pgfmathresult)={\numR}
		
		\forloop[1]{itR}{0}{\value{itR} < \j }
		{
			\pgfmathparse{ int(\theitR+1) } \let\tS\pgfmathresult  	
			\checkvSave(\tS)							
			\pgfmathsetmacro{\save}{\cachedata}  
		
			\pgfmathparse{int(\j-\theitR+1)} \checkvNumeratorL(\pgfmathresult)
			\pgfmathsetmacro{\tmpL}{\cachedata} 	
		    
			\pgfmathparse{ int(\theitR+2) } \checkvNumeratorR(\pgfmathresult)
			\pgfmathsetmacro{\tmpR}{\cachedata} 	       
	
			\pgfmathparse{ int(\theitR+1) } \checkvN(\pgfmathresult)
			\pgfmathparse{ \cachedata / (\tmpL+\tmpR) } \let\tmp\pgfmathresult
			
			\pgfmathparse{ \save + \tmpR * \tmp } 
			\vN(\tS)={\pgfmathresult}

			\pgfmathparse{ \tmpL * \tmp } \let\tmpsave\pgfmathresult
			\pgfmathparse{ int(\theitR+2) } \let\tS\pgfmathresult      
			\vSave(\tS)={\tmpsave}
		}
		
		\pgfmathparse{ int(\j+1) } \checkvSave(\pgfmathresult )
		\vN(\tS)={\cachedata}
	}

	\pgfmathparse{int( #2+1) } \let\lastIndex\pgfmathresult
	\checkvN(1) \pgfmathsetmacro{\first}{\cachedata}  
	\foreach \i [remember=\a as \lasta (initially \first)] in {2,...,\lastIndex}
	{
		\checkvN(\i) \def\a{\lasta,\cachedata} 
		\ifthenelse{\i=\lastIndex}{ \xdef#5{\a} }{}
	}
    
	\foreach \i in {1,...,\lastIndex}
	{
	    \clrarray{vNumeratorR}(\i)
	    \clrarray{vNumeratorL}(\i)
	    \clrarray{vN}(\i)
	    \clrarray{vSave}(\i)
	}
	\delarray\vN
	\delarray\vNumeratorL
	\delarray\vNumeratorR
	\delarray\vSave
}

\newcounter{countvalues}
\newcounter{getBasis}
\newcommand{\bsplinebasis}[4] 
{						
    \newarray\vKnots 	
    \readarray{vKnots}{#1}

    \pgfmathparse{#3}  \let\i\pgfmathresult
    \pgfmathparse{#2}  \let\p\pgfmathresult

    \setcounter{countvalues}{0}
    \setcounter{getBasis}{\p} 
    \pgfmathparse{int(\i+\p)}
    \foreach \knotspan in {\i,...,\pgfmathresult}
    {
	\pgfmathparse{ int(\knotspan+1+1) } \checkvKnots(\pgfmathresult)
	\pgfmathsetmacro{\tmpR}{\cachedata} 						
	
	\pgfmathparse{ int(\knotspan+1) } \checkvKnots(\pgfmathresult) 
	\pgfmathsetmacro{\tmpL}{\cachedata} 						

	\pgfmathparse{ \tmpR - \tmpL } \let\spansize\pgfmathresult  			
   
        \pgfmathparse{ \spansize > 0.0 } \let\bNonZero\pgfmathresult
        \ifthenelse{ \bNonZero = 1 }
        {
            \foreach \percentU in {0,10,...,100}
            {
		\pgfmathparse{\tmpL+\spansize*\percentU/100} \let\u\pgfmathresult

		\bsplinevalue{#1}{\p}{\u}{\knotspan}{\Basis} 
		\def\basisfuncarray{{\Basis}} 				
		
		\pgfmathparse{\basisfuncarray[\thegetBasis]} 
		
		\ifthenelse{\thecountvalues=0}
		{ 
			\xdef\nodeB{"\u,\pgfmathresult"}
		}{
			\xdef\nodeB{\nodeB,"\u,\pgfmathresult"}  
		}
		\addtocounter{countvalues}{1}
            }
        }{}
        \addtocounter{getBasis}{-1}
    }
    
	\xdef#4{\nodeB}

	\delarray\vKnots
}

\newcommand{\bsplinebasisspan}[6] 
{							
    \newarray\vKnots 	
    \readarray{vKnots}{#1}

    \pgfmathparse{#5}  \let\plotknotspan\pgfmathresult
    \pgfmathparse{#4}  \let\splineknotspan\pgfmathresult
    \pgfmathparse{#3}  \let\i\pgfmathresult
    \pgfmathparse{#2}  \let\p\pgfmathresult

    \setcounter{countvalues}{0}
    \setcounter{getBasis}{\i} 
    \foreach \knotspan in {\plotknotspan}
    {
	\pgfmathparse{ int(\knotspan+1+1) } \checkvKnots(\pgfmathresult)
	\pgfmathsetmacro{\tmpR}{\cachedata} 						
	
	\pgfmathparse{ int(\knotspan+1) } \checkvKnots(\pgfmathresult) 
	\pgfmathsetmacro{\tmpL}{\cachedata} 						

	\pgfmathparse{ \tmpR - \tmpL } \let\spansize\pgfmathresult  			
   
        \pgfmathparse{ \spansize > 0.0 } \let\bNonZero\pgfmathresult
        \ifthenelse{ \bNonZero = 1 }
        {
            \foreach \percentU in {0,10,...,100}
            {
		\pgfmathparse{\tmpL+\spansize*\percentU/100} \let\u\pgfmathresult

		\bsplinevalue{#1}{\p}{\u}{\splineknotspan}{\Basis} 
		\def\basisfuncarray{{\Basis}} 				
		
		\pgfmathparse{\basisfuncarray[\thegetBasis]} 
		
		\ifthenelse{\thecountvalues=0}
		{ 
			\xdef\nodeB{"\u,\pgfmathresult"}
		}{
			\xdef\nodeB{\nodeB,"\u,\pgfmathresult"}  
		}
		\addtocounter{countvalues}{1}
            }
        }{}
        \addtocounter{getBasis}{-1}
    }
    
	\xdef#6{\nodeB}

	\delarray\vKnots
}

\newcommand{\plotbsplinebasis}[4] 	
{							
	\bsplinebasis{#1}{#2}{#3}{\nodeOut}
	\def\nodearray{{\nodeOut}}

	\xdef\name{ }
	\addtocounter{countvalues}{-1}
	\foreach \i in {0,...,\thecountvalues}
	{
		\pgfmathparse{\nodearray[\i]}
		\coordinate (point\i) at (\pgfmathresult);	  
		\xdef\name{ \name (point\i)  }
	}
	
	\draw[#4] plot coordinates{ \name };
	
	\xdef\name{ }
}

\newcommand{\plotbsplinesegment}[6] 	
{								
	\bsplinebasisspan{#1}{#2}{#3}{#4}{#5}{\nodeOut}
	\def\nodearray{{\nodeOut}}

	\xdef\name{ }
	\addtocounter{countvalues}{-1}
	\foreach \i in {0,...,\thecountvalues}
	{
		\pgfmathparse{\nodearray[\i]}
		\coordinate (point\i) at (\pgfmathresult);	  
		\xdef\name{ \name (point\i)  }
	}
	
	\draw[#6] plot coordinates{ \name };
	
	\xdef\name{ }
}

\newcommand{\plotbsplineaccumulated}[5] 	
{						
    \newarray\vKnots 	
    \readarray{vKnots}{#1}
    \newarray\vSubCoef 	
    \readarray{vSubCoef}{#3}
    
    \pgfmathparse{#4}  \let\plotknotspan\pgfmathresult
    \pgfmathparse{#4}  \let\splineknotspan\pgfmathresult
    \pgfmathparse{#2}  \let\p\pgfmathresult
    
    \setcounter{countvalues}{0}
    \pgfmathparse{int( \p+1) } \let\lastIndex\pgfmathresult
    \foreach \knotspan in {\plotknotspan}
    {
        \pgfmathparse{ int(\knotspan+1+1) } \checkvKnots(\pgfmathresult)
        \pgfmathsetmacro{\tmpR}{\cachedata} 						
        
        \pgfmathparse{ int(\knotspan+1) } \checkvKnots(\pgfmathresult) 
        \pgfmathsetmacro{\tmpL}{\cachedata} 						
        
        \pgfmathparse{ \tmpR - \tmpL } \let\spansize\pgfmathresult  			
        
        \pgfmathparse{ \spansize > 0.0 } \let\bNonZero\pgfmathresult
        \ifthenelse{ \bNonZero = 1 }
        {
            \foreach \percentU in {0,10,...,100}
            {
                \pgfmathparse{\tmpL+\spansize*\percentU/100} \let\u\pgfmathresult
                
                \bsplinevalue{#1}{\p}{\u}{\splineknotspan}{\Basis} 
                \def\basisfuncarray{{\Basis}} 				
                
                \setcounter{getBasis}{0} 
                \pgfmathparse{\basisfuncarray[\thegetBasis]} 
                \let\basisValue\pgfmathresult
                
                \checkvSubCoef(1) \pgfmathsetmacro{\coef}{\cachedata}  
                \pgfmathparse{ \basisValue * \coef } \let\first\pgfmathresult
                
                \xdef\lastx{\first}
                \foreach \i in {2,...,\lastIndex}
                { 
                    \addtocounter{getBasis}{1}       
                    \pgfmathparse{\basisfuncarray[\thegetBasis]}
                    \let\basisValue\pgfmathresult
                    
                    \checkvSubCoef(\i) \pgfmathsetmacro{\coef}{\cachedata}  
                    \pgfmathparse{ \lastx + \basisValue * \coef } \let\sum\pgfmathresult
                    
                    \xdef\lastx{\sum}
                    
                    \ifthenelse{\i=\lastIndex}
                    {
                        \ifthenelse{\thecountvalues=0}
                        { 
                            \xdef\nodeBB{"\u,\sum"}
                        }{
                            \xdef\nodeBB{\nodeBB,"\u,\sum"}  
                        }
                        \addtocounter{countvalues}{1}
                    }{}
                    
                }
            }
        }{}
    }
    
    \delarray\vKnots
    \delarray\vSubCoef
    
    \def\nodearray{{\nodeBB}}
    
    \xdef\name{ }
    \addtocounter{countvalues}{-1}
    \foreach \i in {0,...,\thecountvalues}
    {
        \pgfmathparse{\nodearray[\i]}
        \coordinate (point\i) at (\pgfmathresult);                
        \xdef\name{ \name (point\i)  }
    }
    
    \draw[#5] plot coordinates{ \name };
    
    \xdef\name{ }
}

\begin{document}
    
\title{Algorithms for geometrical operations with NURBS surfaces.}
\begin{frontmatter}

\author[ifbaddr]{Gernot Beer\corref{cor1}}

\address[ifbaddr]{Institute of Structural Analysis, Graz University
  of Technology, Lessingstraße 25/II, 8010 Graz, Austria}


\cortext[cor1]{Corresponding author.
  mail: \url{gernot.beer@tugraz.at}, web: \url{www.ifb.tugraz.at}}

\begin{abstract}
The aim of the paper is to show algorithms for geometrical manipulations on NURBS surfaces. These include generating NURBS surfaces that pass through given points, calculating the minimum distance to a point and include line to surface and surface to surface intersections.
\end{abstract}

\end{frontmatter}

\section{Introduction}
Nonuniform rational B-splines or NURBS have been used by the Computer Aided Design (CAD) community for decades. The reason for this is that they are very suitable for defining geometrical shapes and for geometrical operations. A great number of publications on NURBS exist, here we quote the NURBS book \cite{Piegl1997b} and a paper \cite{DIMAS1999741} that give a good summary. With the publication of the book Isogeometric Analysis \cite{Cottrell2009} it was suggested that NURBS would also benefit numerical simulation methods such as the Finite Element (FEM)and Boundary Element method (BEM). This was followed by a number of papers discussing the implementation, for example (\cite{Hughes2005a,Bazilevs2010a,Borden2010a,simpson2012two,Trevelyan2013,scott2013isogeometric}). A book on the implementation of isogeometric methods into BEM simulation programs was published \cite{BeerMarussig}.

A proper implementation of NURBS into FEM and BEM programs, however required revisiting of some geometrical operations. While the CAD community is mainly interested in the graphical display of geometry, the simulation community is interested in generating an 'analysis suitable' geometry, i.e. one that allows a suitable volume or boundary discretisation.
In this paper we concentrate of surfaces, as they are used in BEM, and present algorithms for some geometrical operations.
It should be noted that CAD programs have very sophisticated algorithms for computing surface to surface intersections. However, these algorithms produce data that are not 'analysis suitable'.

\section{B-splines and NURBS}

 B-splines are an attractive alternative to Lagrange polynomials and Serendipity functions predominantly used in simulation. The basis for creating the functions is the \textit{knot vector}. This is a vector containing a series of non-decreasing values of the local coordinate $\xi$:
\begin{equation}
\label{ }
\Xi=\left(\begin{array}{cccc}\xi_{0} & \xi_{1} & \cdots & \xi_{N}\end{array}\right)
\end{equation}

A B-spline basis function of order $p=0$ (constant) is given by:
\begin{equation}
B_{i,0}(\mathrm{\xi}) =	\left\{ \begin{array}{l}
			1 \quad \text{if}  \quad  \mathrm{\xi}_{i}\leqslant \mathrm{\xi} < \xi_{i+1} \\
			0 \quad \text{otherwise} \\
			\end{array} 
		\right.
\label{Bf1}
\end{equation}

Higher order basis functions are defined by referencing lower order functions:
\begin{equation}
  B_{i,p}(\mathrm{\xi})=\frac{\xi-\xi_{i}}{\xi_{i+p}-\xi_{i}} \: B_{i,p-1}(\xi) 
			   + \frac{\xi_{i+p+1}-\xi}{\xi_{i+p+1}-\xi_{i+1}} \: B_{i+1,p-1}(\xi)
\label{Bf2}
\end{equation}
B-spines are associate with \textit{anchors}.
The location of the $i$-th anchor in the parameter space can be computed by:
\begin{equation}
\label{Greville}
\abscissa_i = \frac{\uu_{i+1}+\uu_{i+2} + \dots +\uu_{i+p}}{p} \qquad i=0,1, \dots ,I-1.
\end{equation}

\bigskip

Nonuniform rational B-splines or NURBS are based on B-splines but have improved properties for the definition of geometry.
NURBS of order $p$ are defined as:
\begin{equation}
\label{NURBS1}
  N_{i}^p(\xi)=\frac{B_{i,p}(\xi) \: \mathrm{w}_{i}}{\sum_{j=0}^{I}B_{j,p}(\xi) \: \mathrm{w}_{j}}  
\end{equation}
where $I+1$ is the number of basis functions and $\mathrm{w}_{i}$ are weights. 
This can be extended to two dimensions using a tensor product:
\begin{equation}
\label{NURBS2}
  N_{ij}^{p,q}(\xi,\eta)=  N_{i}^p(\xi) \  N_{j}^q(\eta)
\end{equation}
where $N_{j}^q(\eta)$ is a NURBS of order $q$ in $\eta$ direction.

The coordinates of a point on a curve with the local coordinate $\xi$ can be computed by:
\begin{equation}
\label{fitc}
\mathbf{ x}(\xi)= \sum_{j=1}^{J} R_{j}(\xi) \mathbf{ c}_{j}
\end{equation}
where $R_{j}(\xi_i)$ are NURBS basis functions defined in equation \ref{NURBS1} except that numbering starts at 1 instead of 0. $\mathbf{ c}_{j}$ are control point coordinates.
The vector tangential in direction $\xi$ can be computed by:
\begin{equation}
\label{deriv1}
\mathbf{ v}_1=\frac{\partial \mathbf{ x}(\xi)}{\partial \xi}= \sum_{j=1}^{J} \frac{\partial R_{j}(\xi)}{\partial \xi} \mathbf{ c}_{j}
\end{equation}

The coordinates of a point on a surface with the local coordinate $\xi,\eta$ can be computed by:
\begin{equation}
\label{fits}
\mathbf{ x}(\xi, \eta)= \sum_{j=1}^{J} R_{j}(\xi,\eta) \mathbf{ c}_{j}
\end{equation}
where $R_{j}(\xi,\eta)$ are NURBS basis functions defined in equation \ref{NURBS2}, except that numbering starts at 1 instead of 0. The basis functions are numbered consecutively with a single index $j$ instead of two indices $i$ and $j$.
For a surface the vector in direction $\eta$ can be additionally computed by:
\begin{equation}
\label{deriv2}
\mathbf{ v}_2=\frac{\partial \mathbf{ x}(\xi,\eta)}{\partial \eta}= \sum_{j=1}^{J} \frac{\partial R_{j}(\xi,\eta)}{\partial \eta} \mathbf{ c}_{j}
\end{equation}

\begin{figure}
\begin{center}
\begin{overpic}[scale=0.5]{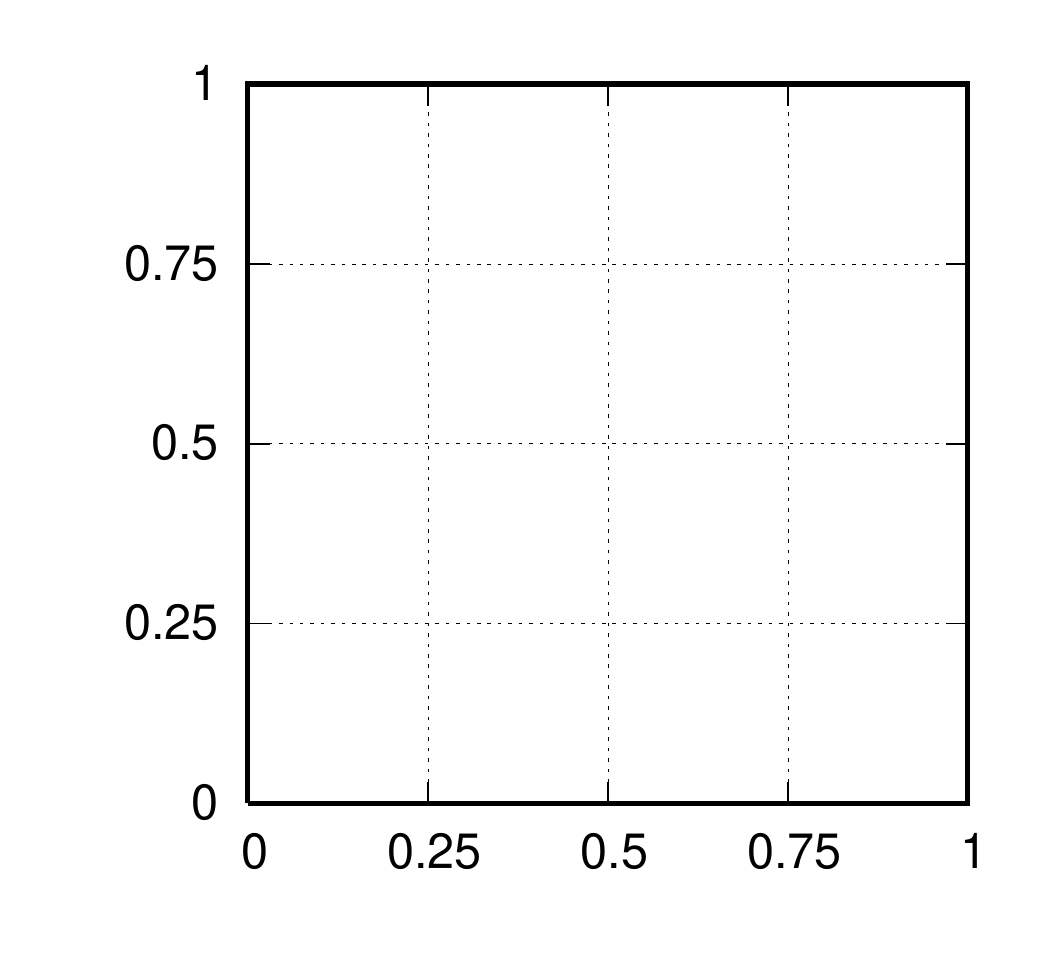}
 \put(50,0){$\xi$}
 \put(0,50){$\eta$}
\end{overpic}
\begin{overpic}[scale=0.5]{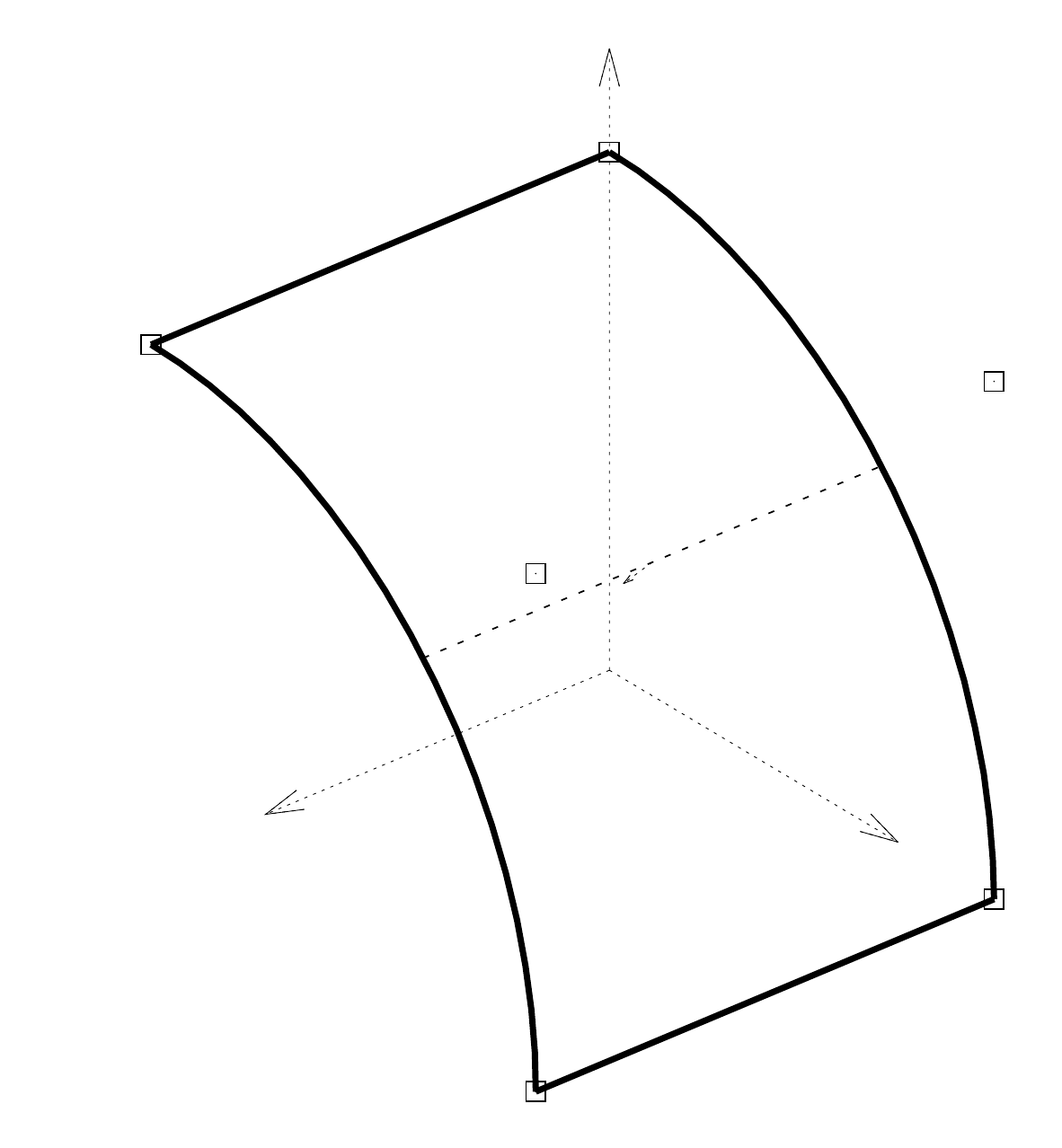}
 \put(80,25){$y$}
 \put(25,32){$x$}
 \put(55,95){$z$}
\end{overpic} 
\caption{A NURBS surface left in local, right in global coordinate system showing control points as hollow squares.}
\label{Surf}
\end{center}
\end{figure}

As an example we show in \myfigref{Surf} a surface created with the knot vectors:
\begin{eqnarray}
\label{ }
\Xi &=& \left[\begin{array}{cccccc}0 & 0 & 0 & 1 & 1 & 1\end{array}\right] \\
\nonumber
H &=& \left[\begin{array}{cccc}0 & 0 & 1 & 1\end{array}\right]
\end{eqnarray}
i.e. quadratic in $\xi$-direction and linear in $\eta$ direction. The control point coordinates and weights are shown in table \ref{tab:Surf}.
\begin{mytable}
  {H}               
  {Control point coordinates and weights for definition of surface.}  
  {tab:Surf}  
  {cccccc}         
  \mytableheader{ $i$ & $j$ & x & y & z & w }  
0 & 0 & 0 & 1 & 0 & 1 \\
1 & 0 & 0 & 1 & 1 & 0.707 \\
2 & 0 & 0 & 0 & 1 & 1 \\
0 & 1 & 1 & 1 & 0 & 1 \\
1 & 1 & 1 & 1 & 1 & 0.707 \\
2 & 1 & 1 & 0 & 1 & 1 \\
\end{mytable}%

\newpage

\subsection{Trimming}
\begin{figure}
\begin{center}
\begin{overpic}[scale=0.5]{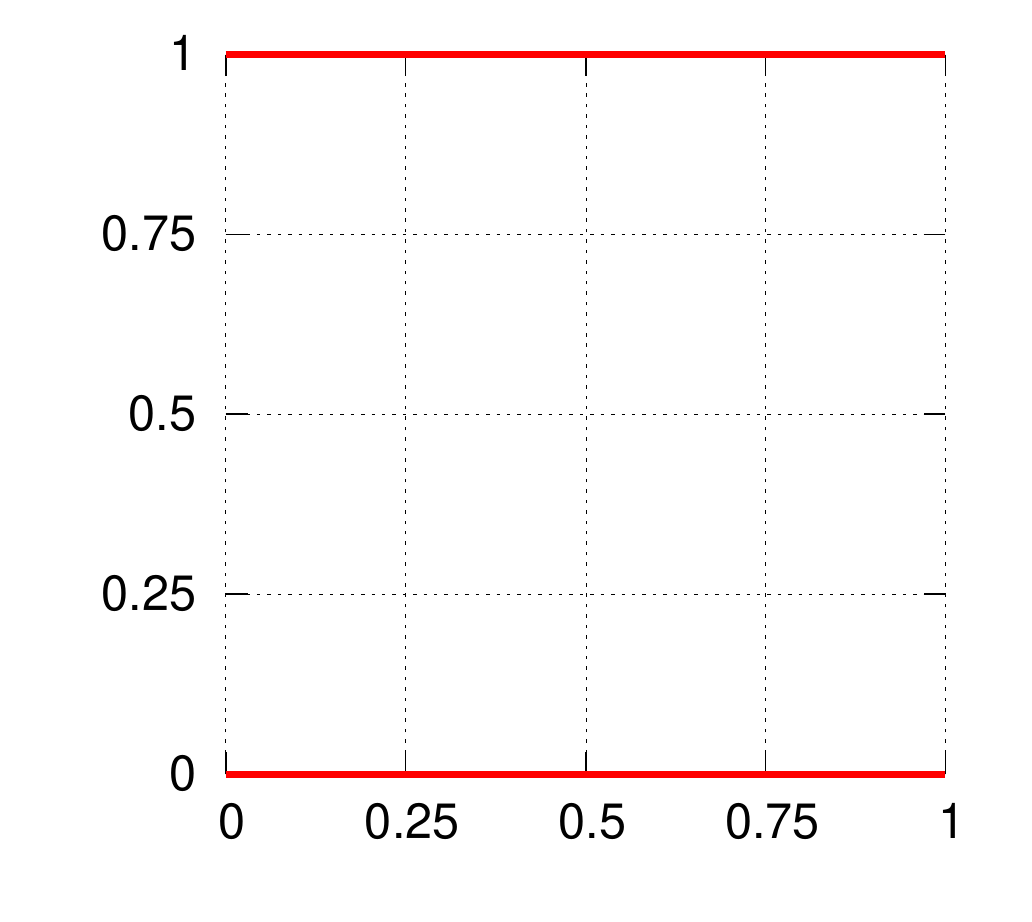}
 \put(50,0){$\xi$}
 \put(0,50){$\eta$}
\end{overpic}
\begin{overpic}[scale=0.5]{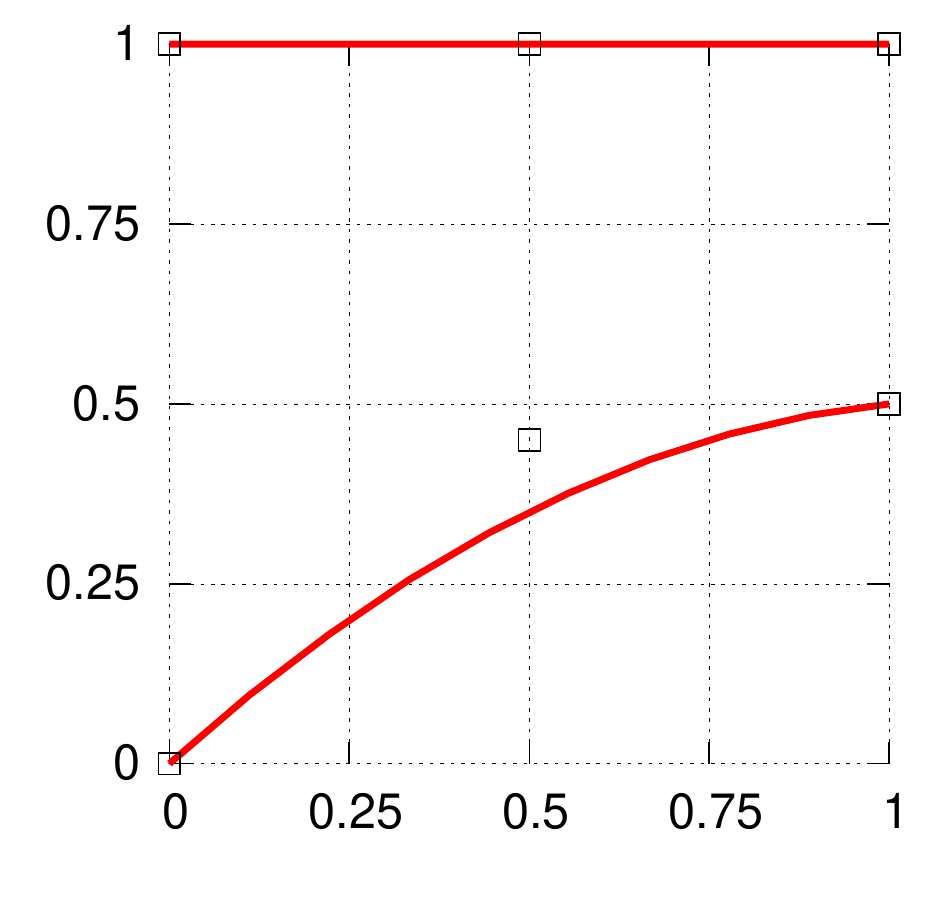}
 \put(50,0){$\hat{\xi}$}
 \put(2,50){$\hat{\eta}$}
\end{overpic}
\begin{overpic}[scale=0.5]{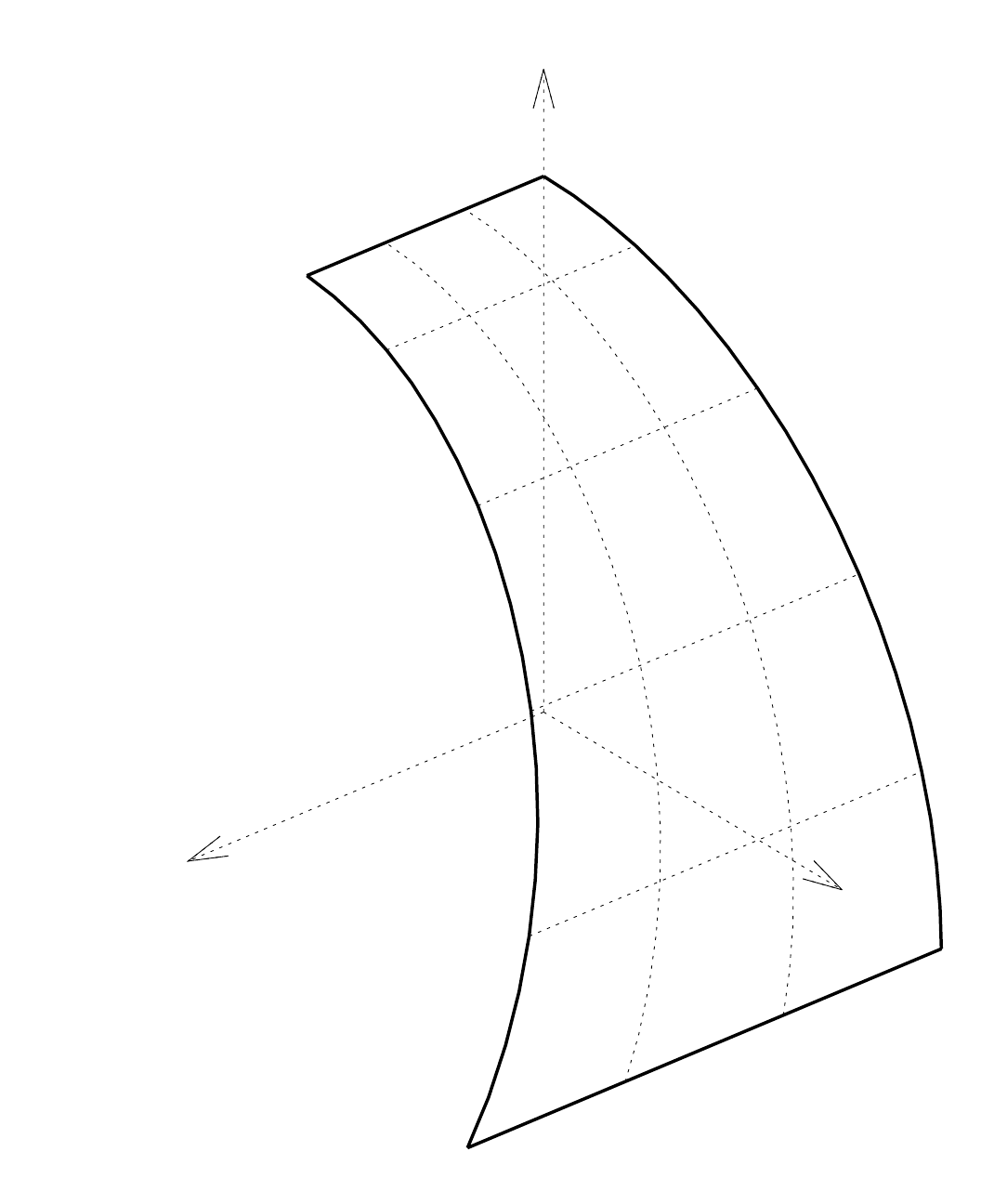}
 \put(75,25){$y$}
 \put(10,25){$x$}
 \put(50,95){$z$}
\end{overpic} 
\caption{Trimming of a surface: 1. Local map , 2. Trimmed map showing trimming curves and control points, 3. Global map of trimmed surface}
\label{SurfT}
\end{center}
\end{figure}
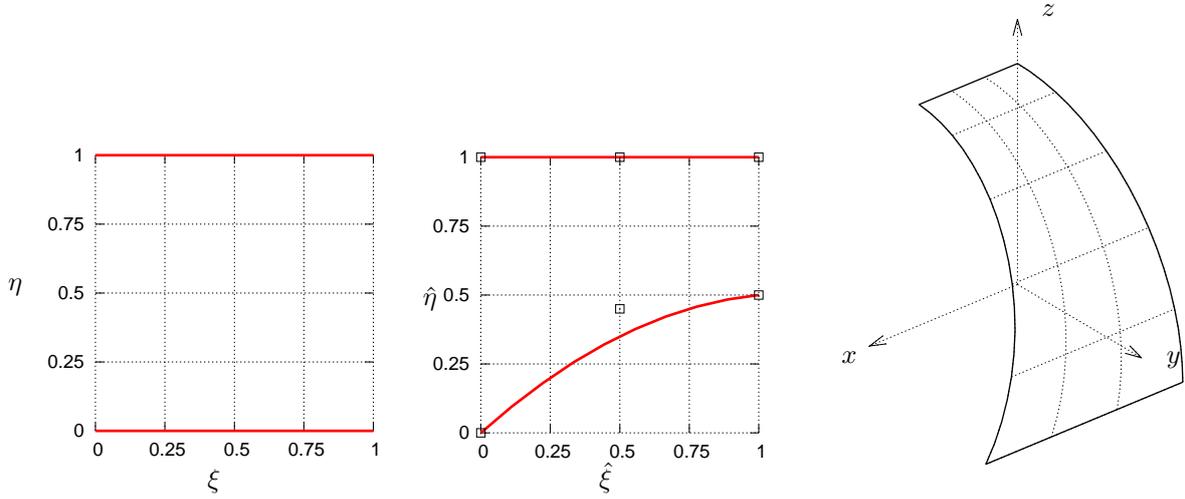
We can cut off a portion of the surface by a process called trimming.
We first define the limits of the trimmed space using basis functions $\hat{R}_{i,p}(\xi) $ and  $ \hat{R}_{j,q}(\eta)$ and control points $ \hat{\textbf{c}}_{ij}$.
We then map from a $\xi,\eta$ space to a local $\hat{\xi},\hat{\eta}$ space:
\begin{equation}
\left\{\begin{array}{c}\hat{\xi} \\ \hat{\eta} \end{array}\right\} =   \sum_{j=0} ^{J} \hat{R}_{j}(\xi,\eta)  \: \hat{\textbf{c}}_{ij}
\label{}
 \end{equation}
Finally we map form the $\hat{\xi},\hat{\eta}$ coordinate system to the $x,y,z$ coordinate system using \myeqref{fits}.
 In \myfigref{SurfT} we show an example on how the surface shown in \myfigref{Surf} can be trimmed.

\newpage

\section{Generating NURBS curves and surfaces through a set of points}
We can determine the NURBS parameters $\mathbf{ c}_{j}$ in such a way that the defined curve or surface either exactly goes through a number of points with the coordinates $\mathbf{ x}(\xi_i , (\eta_{i}))$ or is a least square approximation.
\begin{figure}
\begin{center}
\includegraphics[scale=0.4]{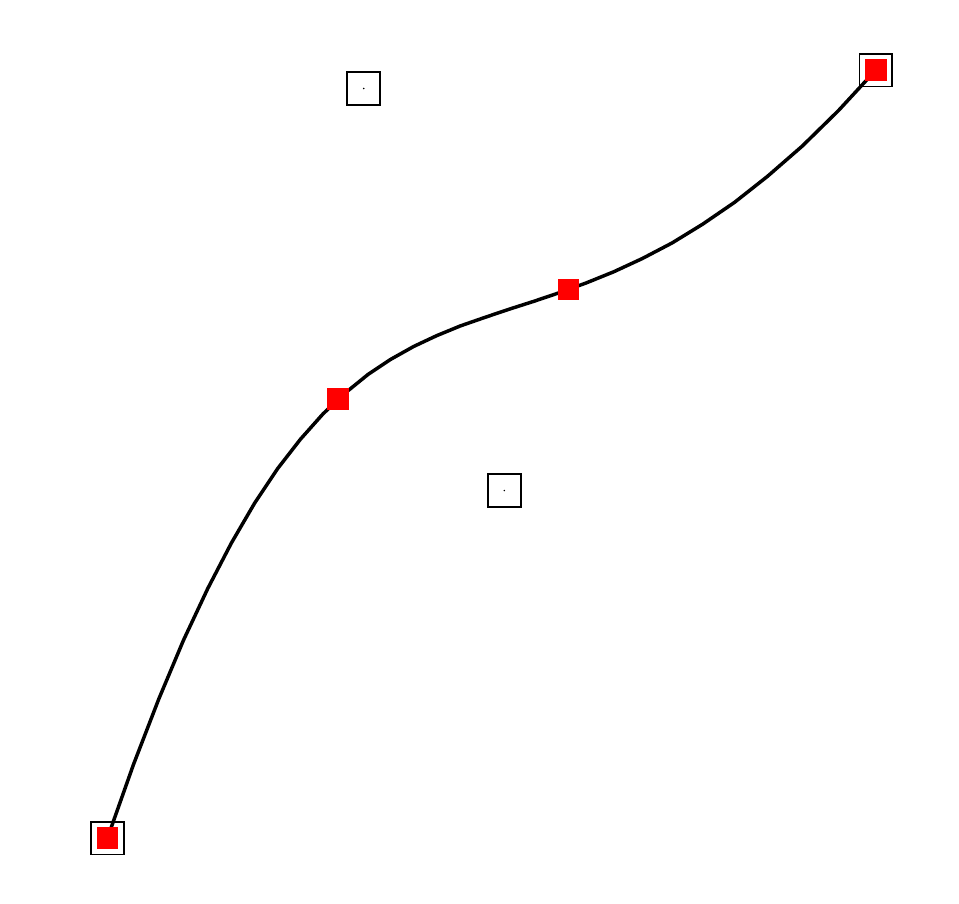}
\includegraphics[scale=0.4]{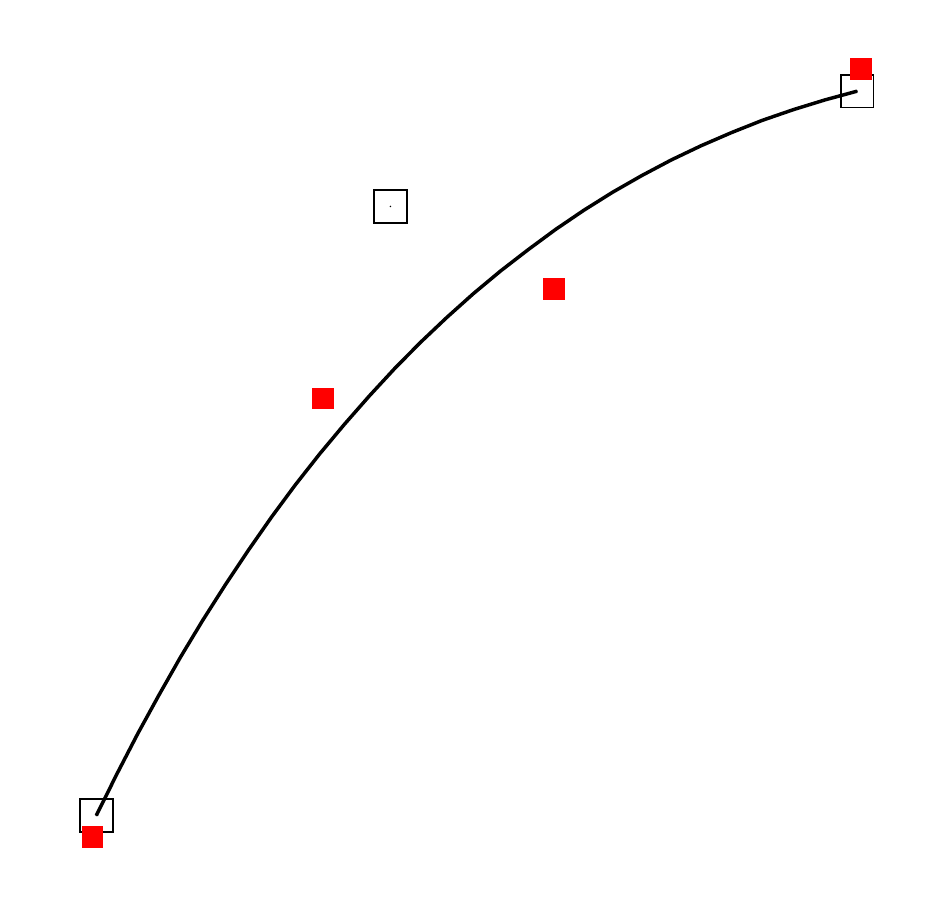}
\caption{Left: curve through 4 points(red squares) using a NURB of order 3, Right: curve through 4 points approximated by a NURB of order 2. Control points are depicted by hollow squares.}
\label{Curvefit}
\end{center}
\end{figure}
\begin{figure}
\begin{center}
\includegraphics[scale=0.5]{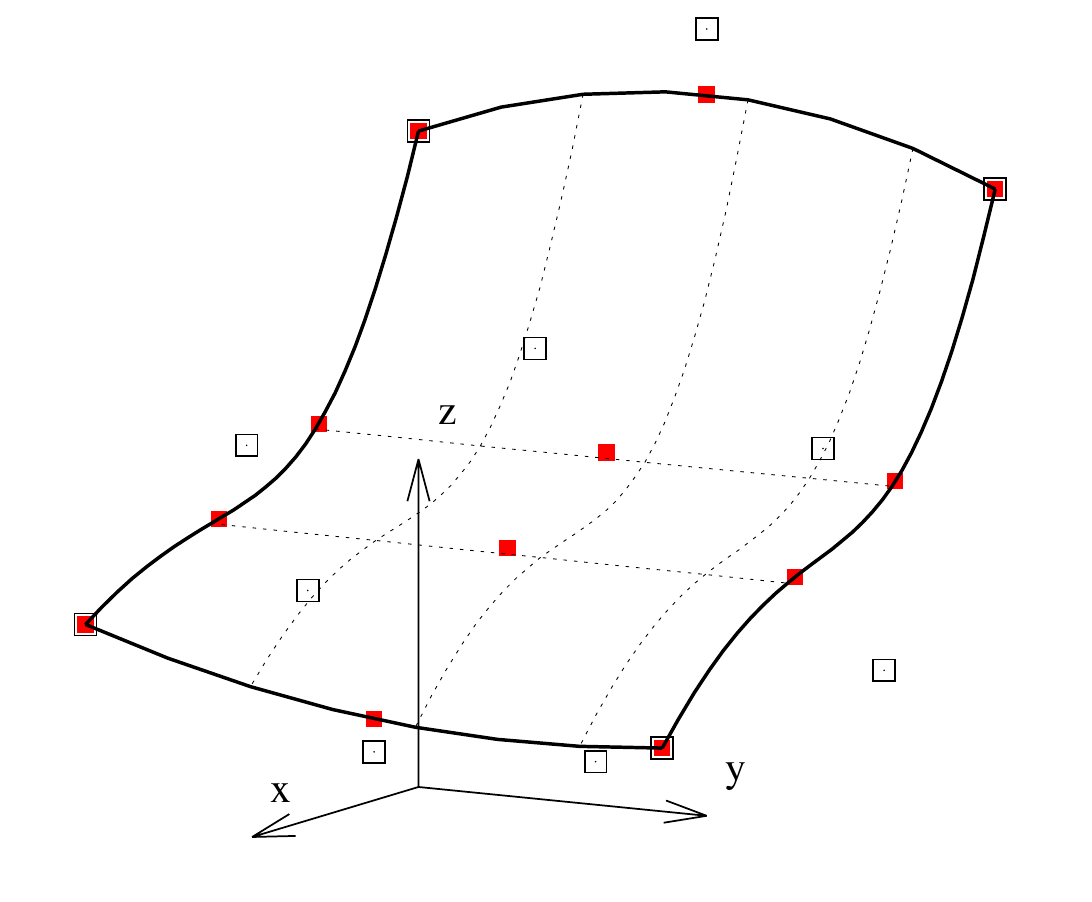}
\includegraphics[scale=0.5]{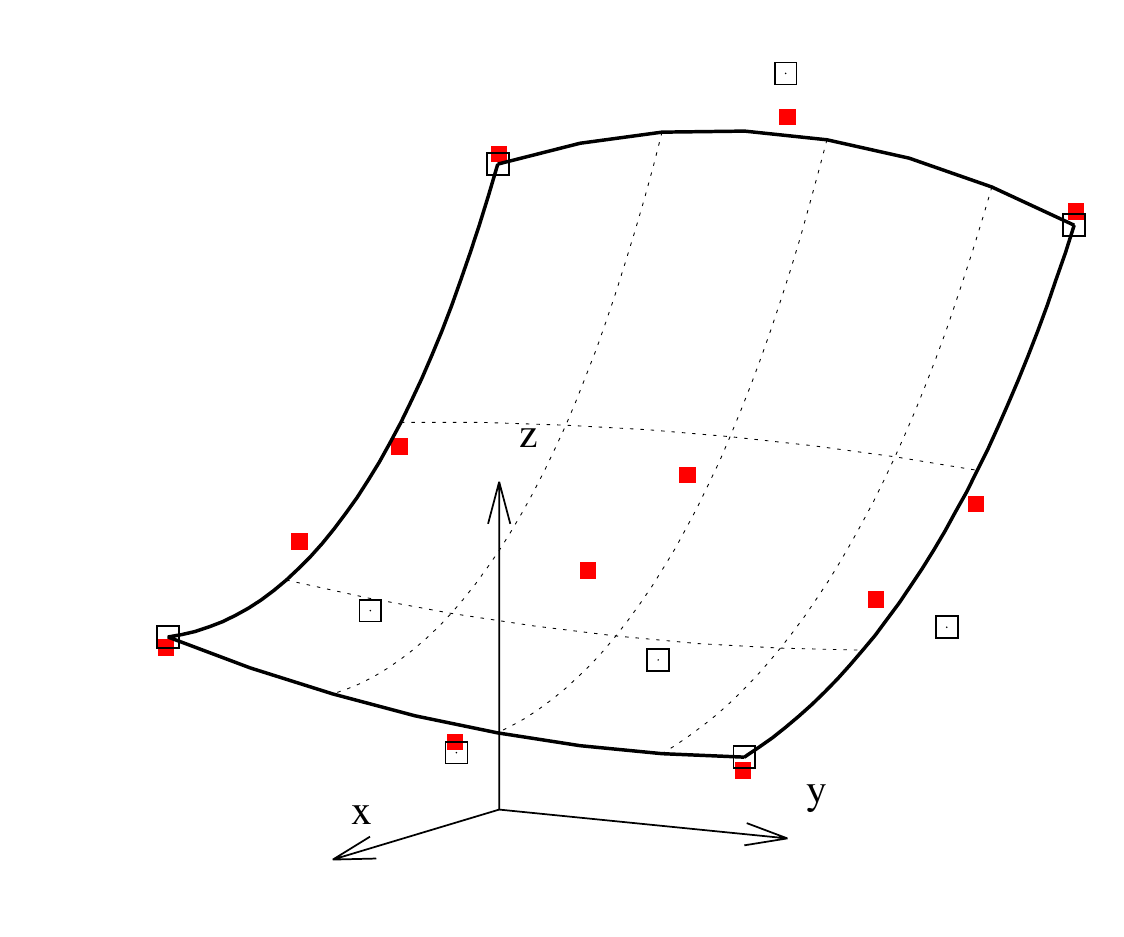}
\caption{Left: surface through 12 points using a NURB of order q=3, p=2, Right: Surface approximating points using NURBS of order q=2,p=2}
\label{CurvefitS}
\end{center}
\end{figure}

First we gather all coordinates of the given points in vector $\{\myVec{x}\}$ and of the unknown control point coordinates in vector $\{\myVec{c}\}$:
\begin{equation}
\label{ }
\{\mathbf{ x}\}= \left\{\begin{array}{c}\mathbf{x}(\xi_{1},(\eta_{1}) )\\ \mathbf{x}(\xi_{2},(\eta_{2}) )  \\\vdots\end{array}\right\} \text{ , } \{\mathbf{c}\}= \left\{\begin{array}{c}\mathbf{c}_{1} \\\mathbf{c}_{2} \\\vdots\end{array}\right\}
\end{equation}

Using (\ref{fitc}) or (\ref{fits}) we arrive at the following matrix equation:
\begin{equation}
\label{ }
 \{\myVec{x}\} = [\mathbf{ A}]\{\myVec{c}\}
\end{equation}

For a curve the matrix $ [\mathbf{ A}]$ is given for example by:
\begin{equation}
\label{ }
 [\mathbf{ A}]= \left(\begin{array}{ccccccc}R_{1}(\xi_{1})& 0 & 0 & R_{2}(\xi_{1})& 0 & 0 & \cdots\\0 & R_{1}(\xi_{1}) & 0 & 0 & R_{2}(\xi_{1}) & 0 & \cdots\\0 & 0 & R_{1}(\xi_{1}) & 0 & 0 & R_{2}(\xi_{1}) & \cdots \\ R_{1}(\xi_{2})& 0 & 0 & R_{2}(\xi_{2})& 0 & 0 & \cdots\\0 & R_{1}(\xi_{2}) & 0 & 0 & R_{2}(\xi_{2}) & 0 & \cdots\\0 & 0 & R_{1}(\xi_{2}) & 0 & 0 & R_{2}(\xi_{2}) & \cdots \\ \vdots& \vdots & \vdots & \vdots& \vdots & \vdots & \ddots\end{array}\right)
\end{equation}

We can solve for the control point coordinates by:
\begin{equation}
\label{ }
\{\myVec{c}\}=  [\mathbf{ A}] \textbackslash \{\myVec{x}\}
\end{equation}
where '$\textbackslash$' denotes either a solution by inverting  $[\mathbf{ A}]$ or a least square approximation of $\{\myVec{x}\}$ depending on the number of points. If the number of points is equal to the number of control points then the curve/surface will go exactly through the specified points.
 \begin{myalgorithm}{Algorithm for computing control point coordinates through a set of points}{alg:aprox}
	\REQUIRE  NURBS functions $R_j$ and local coordinates $\xi_i, (\eta_{i})$, Vector of point coordinates $\{\myVec{x}\}$
	\STATE assemble matrix $[\mathbf{ A}]$
	\STATE solve for control point coordinates $\{\myVec{c}\}=  [\mathbf{ A}] \textbackslash \{\myVec{x}\}$
	\RETURN control point coordinates $\{\myVec{c}\}$
\end{myalgorithm}

Algorithm \ref{alg:aprox} summarises the steps.
An example of curve fitting is shown in \myfigref{Curvefit} of surface fitting in \myfigref{CurvefitS}

\section{Calculating the closest distance between a point and a surface}
To calculate the closest distance of a point to a surface we use an iterative technique.
We start with a point on the surface and compute the vector to the specified point. Using dot products with vectors $\mathbf{ v}_1$ and $\mathbf{ v}_2$ we estimate increments in the local coordinates.
 We iterate until the distance from the specified point to the surface does not change. For surfaces which are very curved we need to limit the increments to a maximum value using a damping coefficient (about 0.5 to 0.75) otherwise the iteration may diverge.

\begin{myalgorithm}{Algorithm for computing minimum distance to a surface}{alg:Mindist}
	\REQUIRE Coordinates of point $\pt{x}_{P}$, NURBS of surface,\\
        Max.~number of iterations $I$ , damping $damp$, Tolerance.
        \STATE Set max. increments: $max_{\xi}=1 , max_{\eta}=1$
	\STATE Select starting point coordinates $\uu_{0},\vv_{0}$ and compute point on surface $\pt{x}_{S0}$
	\STATE Compute $R_{0}= \| \pt{x}_{P} - \pt{x}_{S0}\|$
	\FOR{$i=0$ \TO  $I$}
	\STATE Compute unit vectors $ \myVec{v}_{1}$ and $ \myVec{v}_{2}$ in $\uu$- and $\vv$-directions.
	\STATE $L_\uu= |\myVec{v}_{1}| , L_\vv=|\myVec{v}_{2}|$ 
	\STATE $\triangle \uu = \myVec{v}_{1} \cdot (\pt{x}_P - \pt{x}_{Si}) / L_{\uu}$
        \STATE $\triangle \vv = \myVec{v}_{2} \cdot (\pt{x}_P - \pt{x}_{Si}) / L_{\vv}$
        \STATE $\uu_{i+1}= \uu_{i} + \triangle \uu$
        \IF{$\uu_{i+1}$ >  $max_{\xi}$}
         \STATE $\uu_{i+1}= max_{\xi}$
        \ELSIF{$\uu_{i+1}$ < 0}
         \STATE $\uu_{i+1}=0$
        \ENDIF
        \STATE $\vv_{i+1}= \vv_{i} + \triangle \vv$
         \IF{$\vv_{i+1}$ > $max_{\eta}$}
         \STATE $\vv_{i+1}= max_{\eta}$
        \ELSIF{$\vv_{i+1}$ < 0}
         \STATE $\vv_{i+1}=0$
        \ENDIF
        \STATE Compute new point $\pt{x}_{Si+1}(\uu_{i+1}, \vv_{i+1})$ on surface.
        \STATE Compute distance $R_{i+1} = \| \pt{x}_{P} - \pt{x}_{Si+1}\|$
        \IF{$R_{i+1} - R_{i} < $Tolerance}
        \STATE Exit loop
        \ENDIF
        \STATE Set $max_{\xi}=max_{\xi}*damp$
        \STATE Set $max_{\eta}=max_{\eta}*damp$
	\ENDFOR
	\RETURN Minimum distance $R$, local and global coordinates of point on surface.
\end{myalgorithm}

\begin{figure}[H]
\begin{center}
\begin{overpic}[scale=0.70]{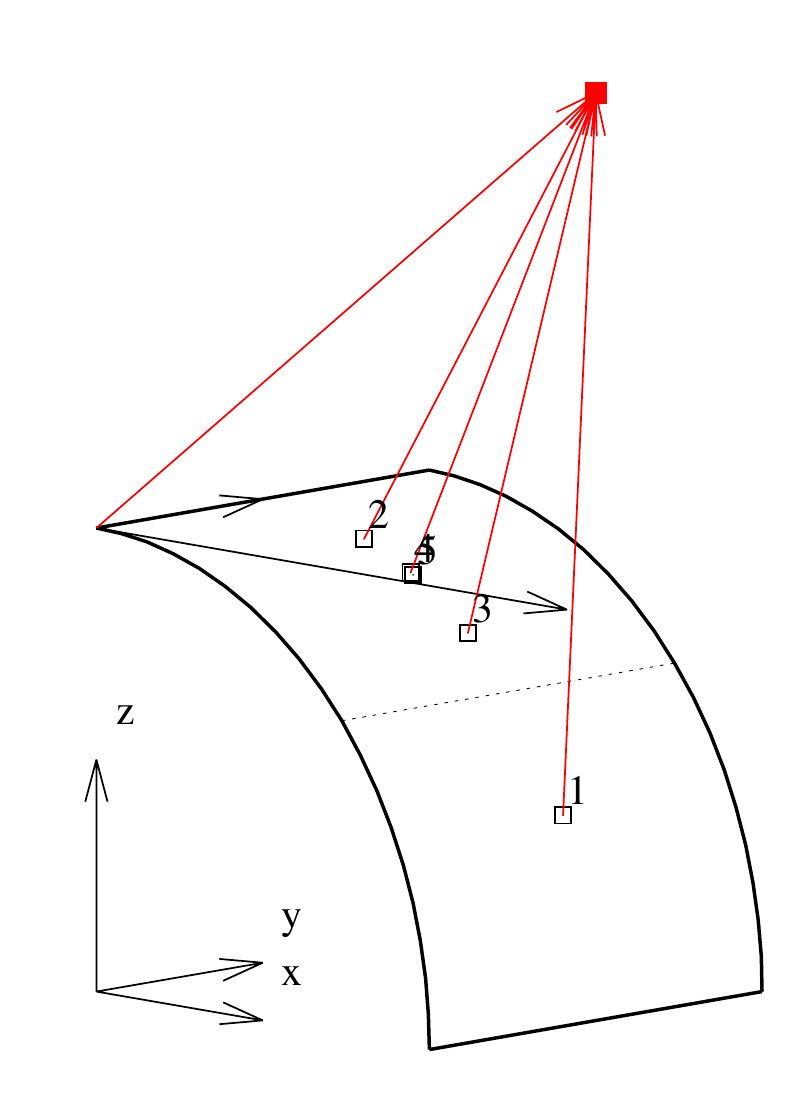}
  \put(50,47){$\mathbf{ v}_1$}
  \put(25,58){$\mathbf{ v}_2$}
\end{overpic} 
\caption{Graphic presentation of the iterative process for the computation of the closest point. Numbers show iteration steps.}
\label{Close}
\end{center}
\end{figure}

Algorithm \ref{alg:Mindist} summarises the steps. \myfigref{Close} illustrates the procedure.

\section{Intersection between a line and a surface}
\begin{figure}[H]
\begin{center}
\begin{overpic}[scale=0.3]{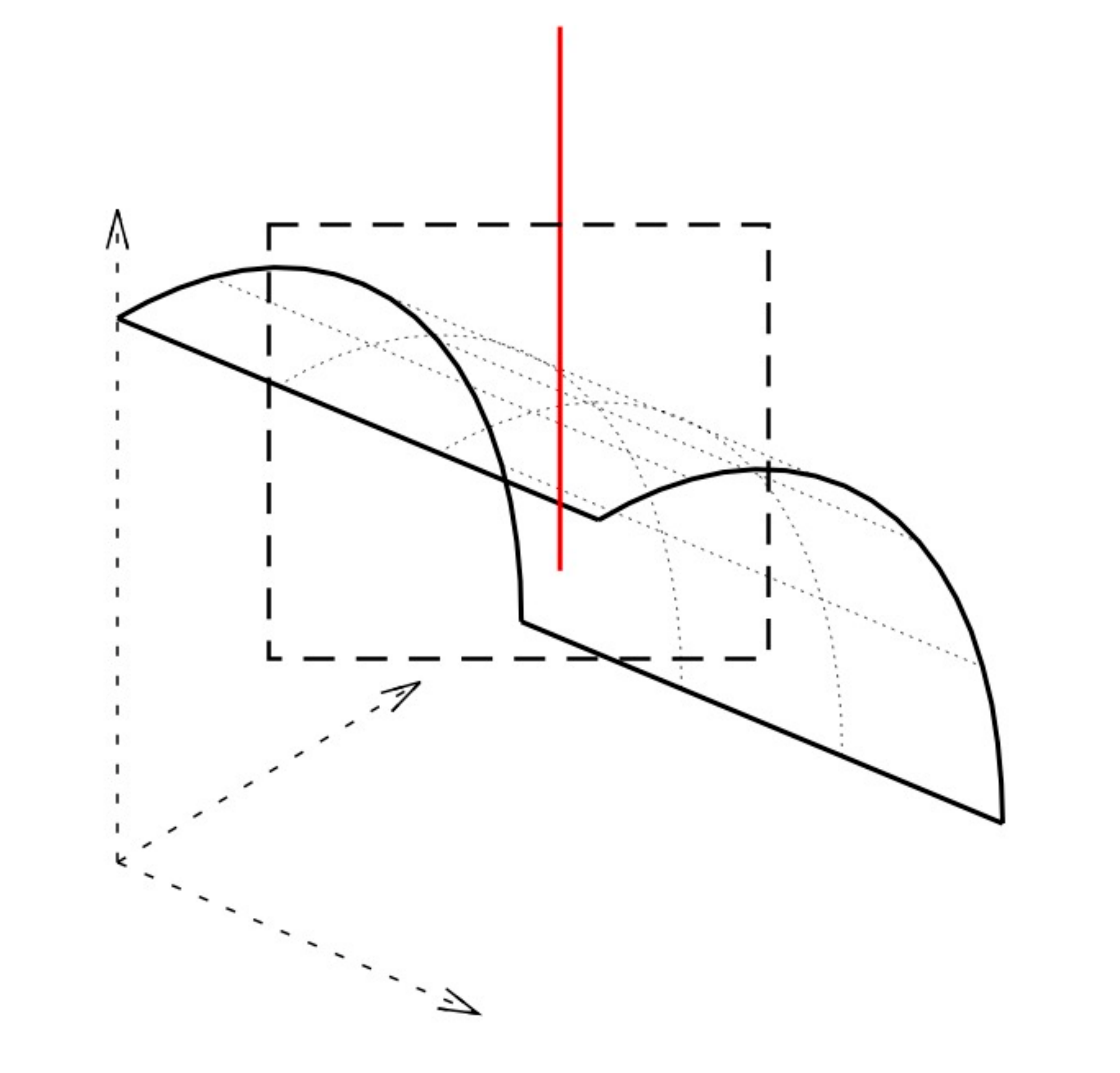}
\end{overpic} 
\begin{overpic}[scale=1.2]{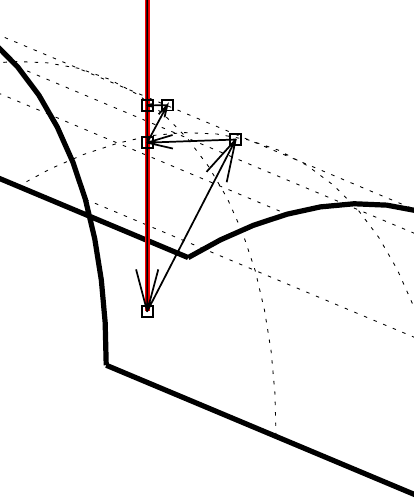}
 \put(33,38){$\mathbf{ x}_0$}
 \put(50,70){$\mathbf{ x}_1$}
 \put(37,78){$\mathbf{ x}_2$}
 \put(20,60){$d\xi_1$}
 \put(20,75){$d\xi_2$}
\end{overpic} 
\caption{Left: A NURBS surface with a line. Right: Detail of iterative procedure to get intersection.}
\label{InsectL}
\end{center}
\end{figure}
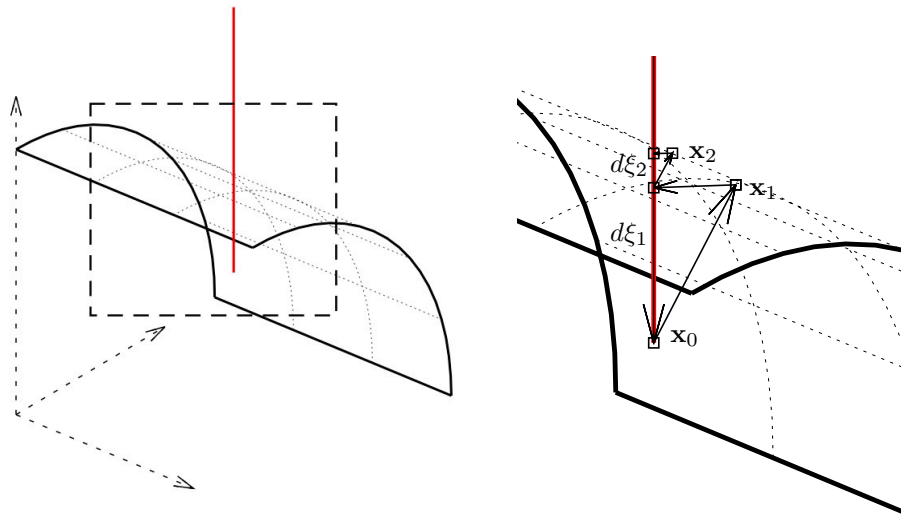

To calculate the intersection of a line with a surface we also use an iterative process. We define the line with a one dimensional NURBS and compute a vector $\mathbf{ V}_{L}$ representing the line. We define $\xi$ as the local coordinate along the line. Starting with $\xi=1$ which corresponds to the end of $\mathbf{ V}_{L}$ (global coordinate $\mathbf{ x}_{0}$) we compute the closest point to the surface as $\mathbf{ x}_{1}$. We then compute and increment of $d\xi_1$ by projecting this point normal to the line,  arriving at a new point on the line closer to the surface. We repeat this process until $d\xi$ is very small, i.e. the point on the line corresponds to a point on the surface.
We explain the procedure in \myfigref{InsectL}.
Algorithm \ref{alg:Insect} outlines the steps.
 \begin{myalgorithm}{Algorithm for computing the intersection of a line with a surface}{alg:Insect}
	\REQUIRE  NURBS of surface, NURBS of line, number of iterations $niter$, tolerance for convergence $tol$
	\STATE compute  vector of line $\mathbf{ V}_{L}$, line length $L$, 
	\STATE Set $\xi=1$  (= end of line).
	\FOR {$i=1$ \TO $niter$}
	 \STATE compute $\mathbf{ x}_{i-1} (\xi)$ on line 
         \STATE  compute $\mathbf{ x}_{i}$ as the point on the surface closest to $\mathbf{ x}_{i-1}$ 
         \STATE  compute $d\xi_{i}= \frac{1}{L} \ (\mathbf{ x}_{i} - \mathbf{ x}_{i-1} ) \cdot  \mathbf{ V}_{L} $
         \STATE compute local coordinate of point on line $\xi= \xi - d\xi_{i}$
         \IF{$d\xi < tol$}
          \STATE Exit 
         \ENDIF
        \ENDFOR
	\RETURN Coordinates of intersection point $\mathbf{ x}_{i}$, local coordinates of point on line $\xi$
\end{myalgorithm}

\newpage

\section{Intersection between a surface and a surface}
\begin{figure}
\begin{center}
\begin{overpic}[scale=0.5]{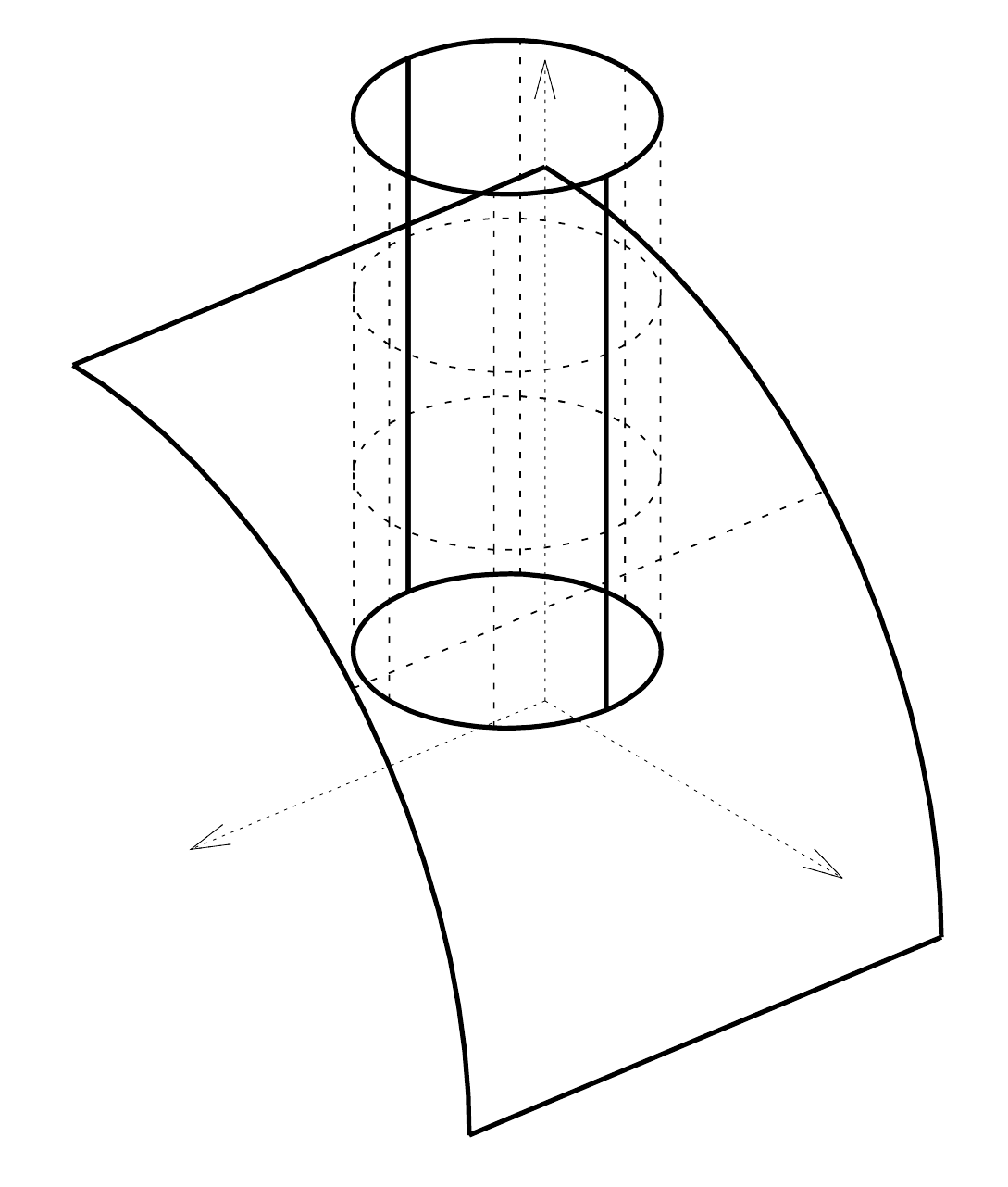}
  \put(10,30){$x$}
  \put(70,30){$y$}
  \put(45,98){$z$}
\end{overpic} 
\caption{Example of 2 intersecting surfaces.}
\label{Surfaces}
\end{center}
\end{figure}

As explained earlier CAD programs have sophisticated and fast intersection algorithms. However, they don't lead to analysis suitable data.
Here we present a simple algorithm that can be used to generate a suitable Boundary element mesh.
We explain the determination of the intersection between 2 surfaces on an example involving the surface defined in section \ref{Surf} and a vertical cylinder depicted in \myfigref{Surfaces}.

\begin{figure}
\begin{center}
\begin{overpic}[scale=0.5]{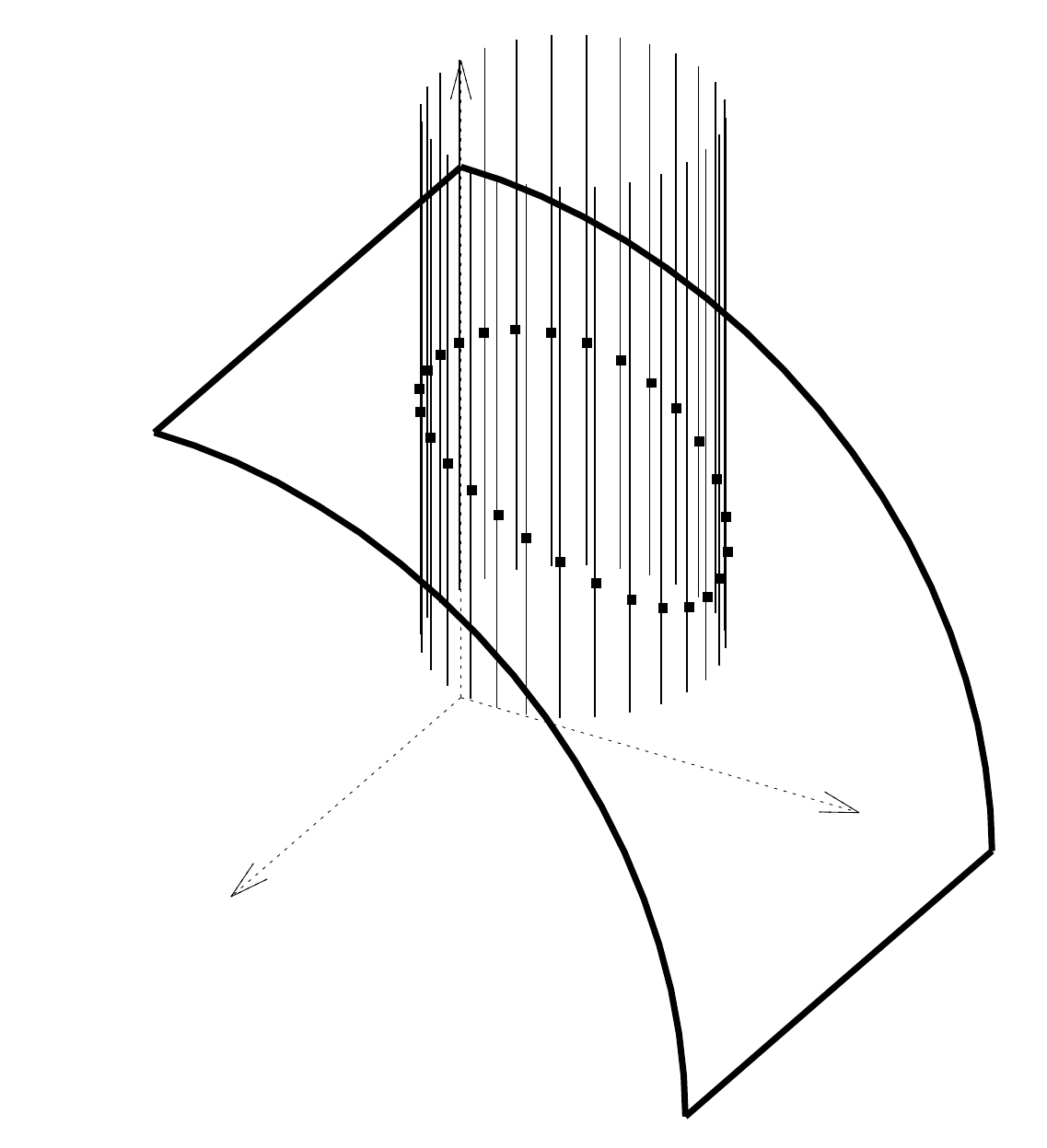}
  \put(15,20){$x$}
  \put(75,30){$y$}
  \put(40,98){$z$}
\end{overpic} 
\begin{overpic}[scale=0.5]{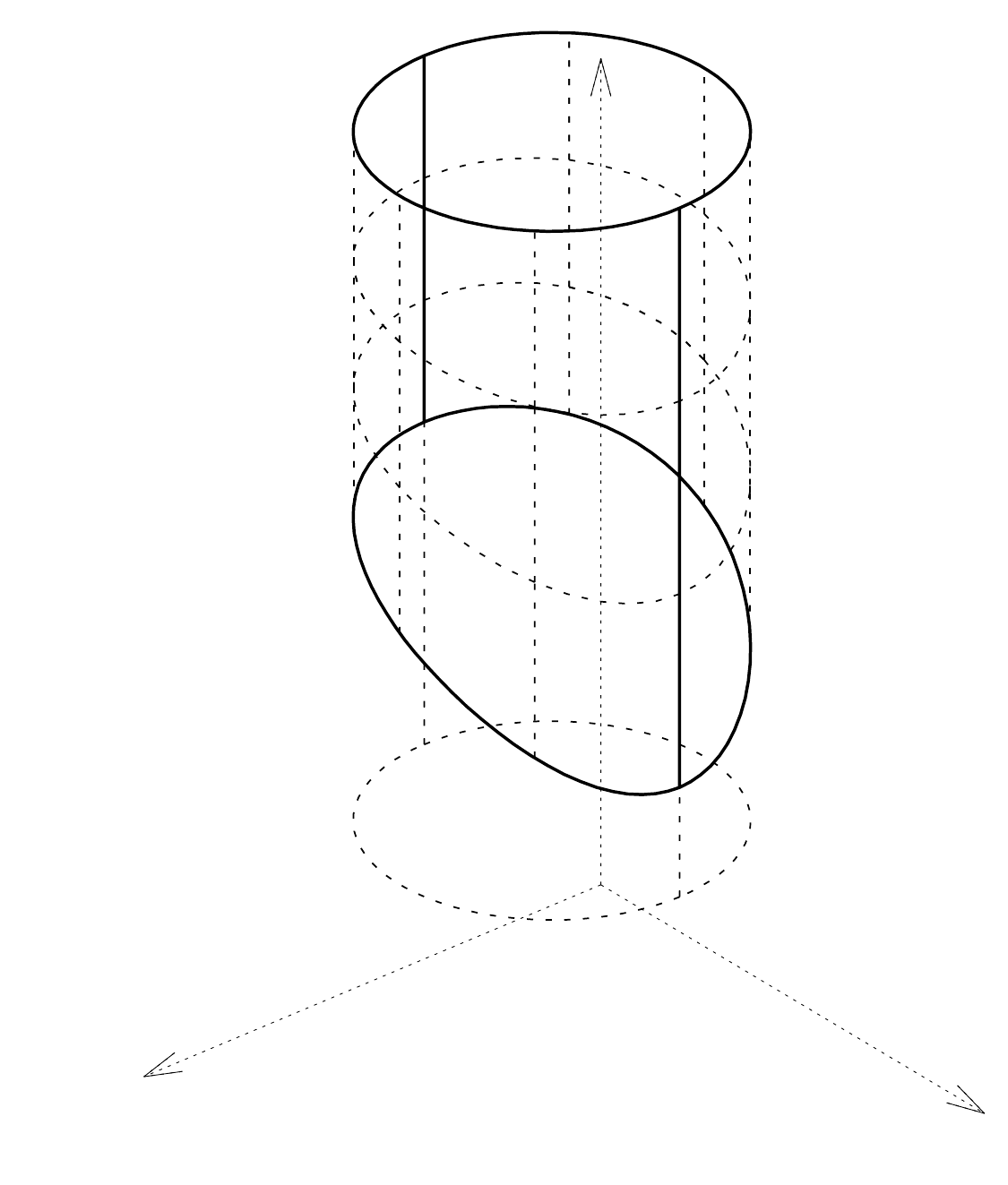}
\end{overpic} 
\caption{Steps 1 and 2 of the process, showing lines and the computed intersection points. Shown on the right is the trimmed intersecting surface}
\label{Step12}
\end{center}
\end{figure}

We specify the cylinder (surface 1) as the \textbf{intersecting} surface and the other surface as the \textbf{intersected} surface (surface 2).
Next we subdivide surface 1 in $\xi$ direction and generate vertical lines as shown in \myfigref{Step12} left. The number of lines will determine the quality of the intersection curve. We use algorithm \ref{alg:Insect} to compute the intersection points with surface 2. Next we use algorithm \ref{alg:aprox} to compute the parameters of a trimming curve that passes through the calculated points. We then trim surface 1 as shown in \myfigref{Step12} right.

\

\remark{These trimming curve parameters can also be taken from data produced by CAD programs.}

\

To produce analysis suitable data the next steps are more involved. We perform these steps in the local $\xi,\eta$ coordinate system of surface 2.
First we divide surface 2 into 4 subsurfaces as shown in \myfigref{Step35} left. For each subsurface we determine trimming curves that go through the computed intersection points and then trim the subsurfaces (\myfigref{Step35} middle). We finally map to global coordinates(\myfigref{Step35} right).
 \begin{figure}
\begin{center}
\begin{overpic}[scale=0.5]{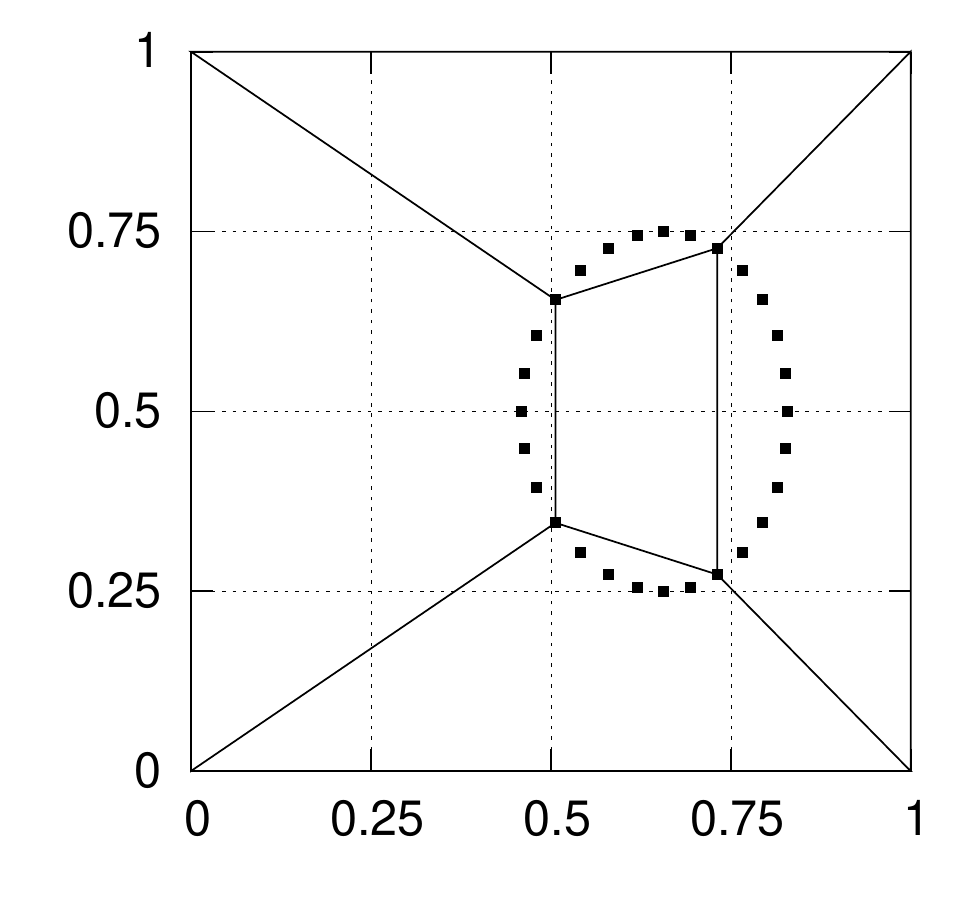}
  \put(50,0){$\xi$}
  \put(0,50){$\eta$}
\end{overpic} 
\begin{overpic}[scale=0.5]{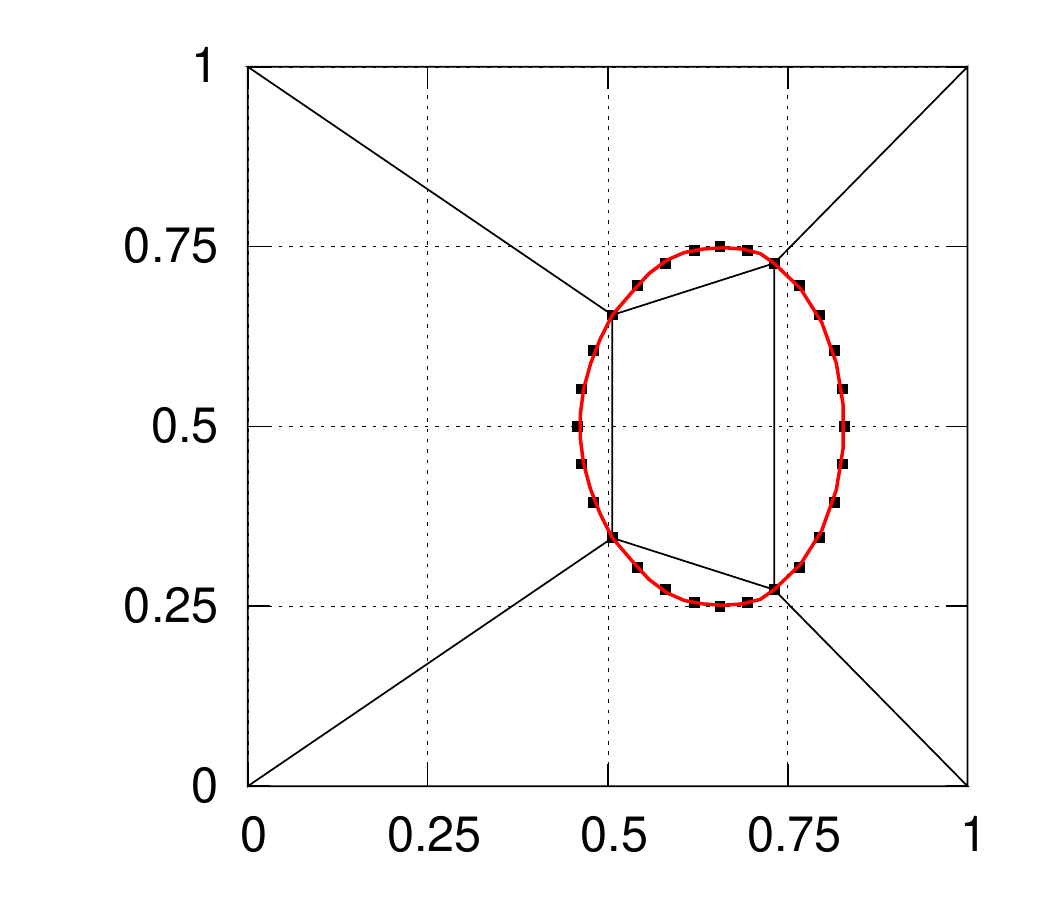}
 \put(50,0){$\xi$}
  \put(0,50){$\eta$}
\end{overpic} 
\begin{overpic}[scale=0.5]{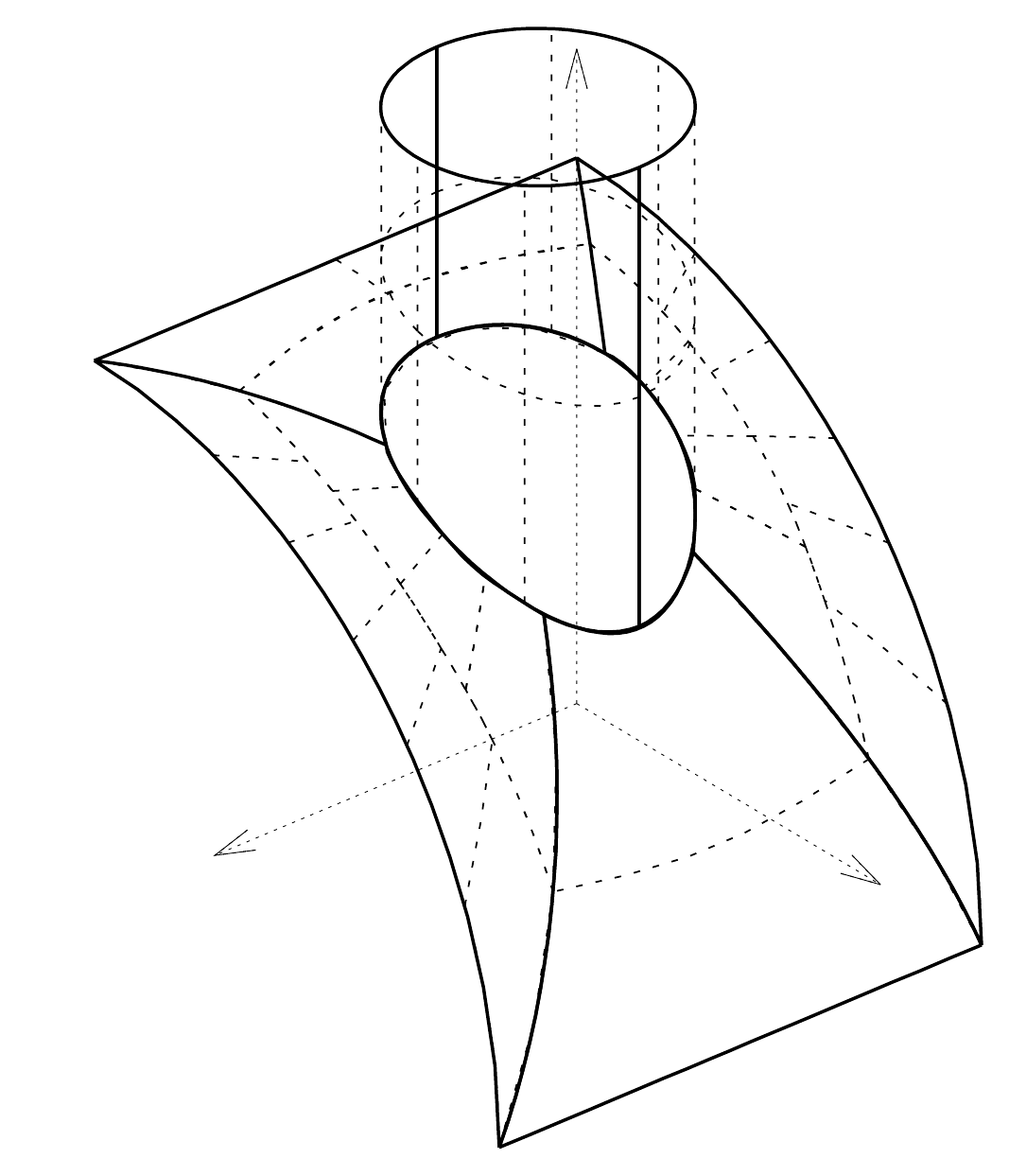}
  \put(15,20){$x$}
  \put(75,30){$y$}
  \put(45,98){$z$}
\end{overpic} 
\caption{Step 3 to 5 of intersection process showing the intersection points in the $\xi,\eta$ coordinate system of the intersected surface.
 The original surface is split into 4 subsurfaces which are trimmed and  then mapped to the global system. The final outcome consisting of 6 surfaces is shown on the right.}
\label{Step35}
\end{center}
\end{figure}
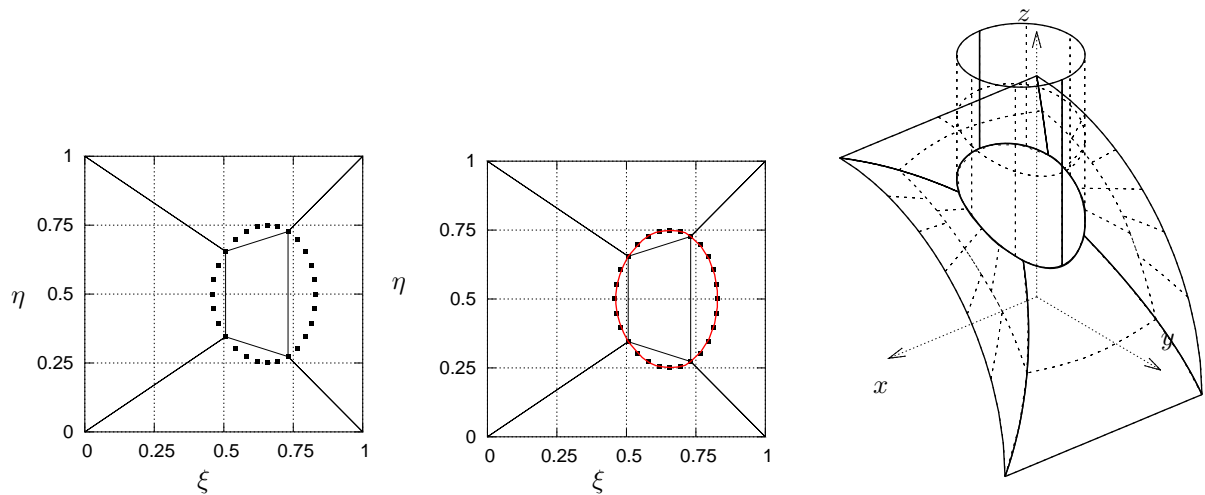
 \begin{myalgorithm}{Algorithm for computing the intersection of a surface with a surface}{alg:InsectS}
	\REQUIRE  NURBS of surface 1 and surface 2, number of intersection lines \textit{nlines}
	\STATE Subdivide surface 1 in $\xi$ direction into \textit{nlines} points.
	\STATE Generate \textit{nlines} lines in $\eta$ direction.
	\STATE Intersect these lines with surface 2 to generate intersection points
	\STATE Use the local $\xi,\eta$ coordinates of intersection points to compute trimming curve for surface 1
	\STATE Trim surface 1
	\STATE Subdivide surface 2 into four subsurfaces in the local coordinate system $\xi,\eta$.
	\STATE Determine trimming curves that go through the intersection points.
	\STATE Trim the subsurfaces
	\STATE Map to global coordinates
	\RETURN  Analysis suitable data
\end{myalgorithm}
Algorithm \ref{alg:InsectS} summarises the steps

\newpage
\section{Summary}
In this short paper some algorithms were presented that the author has used to implement isogeometric  methods into a BEM simulation program. It is hoped that it will be of benefit to developers of simulation software.

\bibliographystyle{myplainnat}
\bibliography{bookbib}

\end{document}